\documentclass[11pt,bezier,amstex,bbm]{article}  

\usepackage{color}              
\usepackage{epsfig}
\usepackage{graphicx}
\usepackage{epsfig}
\usepackage{hyperref}
\usepackage{amsmath,enumerate, amsfonts, amssymb, bbm}
\definecolor{MyDarkBlue}{rgb}{0,0.08,0.50}
\definecolor{BrickRed}{rgb}{0.65,0.08,0}
\newcommand{\rr}[1]{#1} 
\newcommand{\ch}[1]{#1} 

\newtheorem{Lemma}{Lemma}[section]
\newtheorem{Proposition}[Lemma]{Proposition}

\newtheorem{Theorem}[Lemma]{Theorem}

\newtheorem{Corollary}[Lemma]{Corollary}

\newcommand{\proof} {\noindent {\bf Proof}. \hspace{2mm}}

\newcommand{\prob}{\mathbb{P}}

\newcommand{\N}{\mathbb{N}}

\newcommand{\FF}{\mathcal{F}}
\newcommand{\GG}{\mathcal{G}}
\newcommand{\NN}{\mathcal{N}}

\newcommand{\eps}{\varepsilon}

\newcommand{\qed}{\ \ \rule{1ex}{1ex}}

\newcommand{\Zbold}{{\mathbb{Z}}}

\newcommand{\bpi}{\mbox{\boldmath$\pi$}}

\newcommand{\ind}[2]{1_{(e \in \pi(#1,#2))}}

\newcommand{\expec}{\mathbb{E}}
\newcommand{\bfD}{{\bf D}}
\newcommand{\CC}{\mathcal{C}}
\newcommand{\KK}{\mathcal{K}}

\newcommand{\su}{\mathcal{S}}

\newcommand{\KKor}{\mathcal{K}_\infty^{\tt or}}
\newcommand{\KKer}{\mathcal{K}_\infty^{\tt er}}
\newcommand{\KKner}{\mathcal{K}_n^{\tt er}}
\newcommand{\Der}{D^{\tt er}}
\newcommand{\Cer}{C^{\tt er}}
\newcommand{\Ier}{I^{\tt er}}
\newcommand{\Jer}{J^{\tt er}}
\newcommand{\Ner}{N^{\tt er}}

\newcommand{\Ior}{I^{\tt or}}
\newcommand{\Jor}{J^{\tt or}}
\newcommand{\Wor}{W^{\tt or}}
\newcommand{\Wer}{W^{\tt er}}
\newcommand{\GGor}{\GG^{\tt or}}
\newcommand{\GGer}{\GG^{\tt er}}
\newcommand{\Hor}{H^{\tt or}}
\newcommand{\Her}{H^{\tt er}}
\newcommand{\Vor}{V^{\tt or}}
\newcommand{\Ver}{V^{\tt er}}

\newcommand{\UnVer}{A}
\newcommand{\veciota}{\vec{\iota}}
\newcommand\1{\mathbbm{1}}
\renewcommand{\ind}{\1}
\newcommand{\indic}[1]{\1_{\{#1\}}}
\newcommand{\indicwo}[1]{\1_{#1}}
\newcommand{\CBnm}{{\cal B}_{m}^{\sss(n)}}
\newcommand{\Mnm}{M_{m}^{\sss(n)}}

\newcommand{\eqn}[1]{\begin{equation} #1 \end{equation}}
\newcommand{\eqan}[1]{\begin{align} #1 \end{align}}
\newcommand{\lbeq}[1]{\label{#1}}
\newcommand{\refeq}[1]{(\ref{#1})}
\newcommand{\sss}{\scriptscriptstyle}

\newcommand{\op}{o_{\sss \prob}}
\newcommand{\Op}{O_{\sss \prob}}
\newcommand {\convd}{\stackrel{d}{\longrightarrow}}
\newcommand {\convp}{\stackrel{\sss {\mathbb P}}{\longrightarrow}}
\newcommand {\convas}{\stackrel{\sss a.s.}{\longrightarrow}}

\newcommand {\vep}{\varepsilon}
\newcommand{\nn}{\nonumber}

\numberwithin{equation}{section}

\begin{document}
\author{Shankar Bhamidi
\thanks{Department of Mathematics,
The University of British Columbia,
Room 121, 1984 Mathematics Road,
Vancouver, B.C., Canada V6T 1Z2
}
\and
Remco van der Hofstad
\thanks{Department of Mathematics and
Computer Science, Eindhoven University of Technology, P.O.\ Box 513,
5600 MB Eindhoven, The Netherlands. E-mail: {\tt
rhofstad@win.tue.nl}}
\and
Gerard Hooghiemstra
\thanks{DIAM, Delft University of Technology, Mekelweg 4, 2628CD Delft, The
Netherlands, email: g.hooghiemstra@ewi.tudelft.nl} }

\title{Extreme value theory, Poisson-Dirichlet distributions and FPP on random networks
}

\maketitle

\begin{abstract}
We study first passage percolation on the configuration model (CM)
having power-law degrees with exponent $\tau\in [1,2)$. To this
end we equip the edges with exponential weights.  We derive the
distributional limit of the minimal weight of a path between
typical vertices in the network and the number of edges on the
minimal weight path, which can be computed in terms of the
Poisson-Dirichlet distribution. We explicitly describe these
limits via the construction of an infinite limiting object
describing the FPP problem in the densely connected core of the
network. We consider two separate cases, namely, the {\it original
CM}, in which each edge, regardless of its multiplicity, receives
an independent exponential weight, as well as the {\it erased CM},
for which there is an independent exponential weight between any
pair of direct neighbors. While the results are qualitatively
similar, surprisingly the limiting random variables are quite
different.

Our results imply that the flow carrying properties of the network are markedly
different from either the mean-field setting or the locally tree-like setting, which occurs as $\tau>2$, and for which
the hopcount between typical vertices scales as $\log{n}$. In our setting the
hopcount is tight and has an explicit limiting distribution, showing that one can transfer
information remarkably quickly between different vertices in the network.
This efficiency has a down side in that such networks are remarkably fragile to
directed attacks. These results continue a general program by the authors to obtain a
complete picture of how random disorder changes the inherent geometry of various
random network models, see  \cite{SD07,vcg-random-shanky,BhaHofHoo09b}.
\end{abstract}

\vspace{0.6in}


\noindent
{\bf MSC2000 subject classification.}
60C05, 05C80, 90B15.

\section{Introduction}
\label{sec-int}
First passage percolation (FPP) was introduced by Hammersley and Welsh \cite{hamm-welsh}  to
model the flow of fluid through random media. This model has evolved into one of the
fundamental problems studied in modern probability theory, not just for its own sake but
also due to the fact that it plays a crucial role in the analysis of many other
problems in statistical physics, in areas such as the
contact process,  the voter model, electrical resistance problems and in fundamental
stochastic models from evolutionary biology, see e.g. \cite{durrett-lnp}. The
basic model for FPP on (random) graph is
defined as follows: We have some connected graph on $n$ vertices. Each edge is given
some random  weight, assumed to be non-negative,  independent and
identically distributed (i.i.d.) across the edges.
The weight on an edge has the interpretation of the \emph{length} or \emph{cost} of
traversing this edge. Fixing two vertices in the network, we are then interested
in the length and weight of the minimal weight path between these two vertices and the asymptotics of these
statistics as the size of the network tends to infinity.

Most of the classical theorems about FPP deal with the $d$-dimensional
integer lattice, where the connected network is the $[-r,r]^d$ box in the integer lattice and one
is interested in asymptotics of various quantities as $n=(2r+1)^d\to\infty$. In this context,
probabilists are often interested in proving shape theorems, namely, for fixed distance $t$,
showing that $\CC_t/t$ converges to a deterministic limiting set as $t\to \infty$,
where $\CC_t$ is the cluster of all vertices within distance $t$ from the origin.
See e.g., \cite{howard} for a survey of results in this context.

In the modern context such problems have taken on a new significance. The last few years
have witnessed an explosion in the amount of empirical data on networks, including data
transmission networks such as the Internet, biochemical networks such as gene regulatory
networks, spatial flow routing networks such as power transmission networks and
transportation networks such as road and rail networks. This has stimulated an
intense cross-disciplinary effort in formulating network models to understand the
structure and evolution of such real-world networks. Understanding FPP in the context of these random models seems to be of paramount importance,
with the minimal weight between typical vertices representing the cost of transporting
flow between these vertices, while the hopcount, which is defined to be the number of
edges on the minimal weight path between two typical vertices, representing the amount of
time it takes for flow to be transported between these vertices.

In this study we shall analyze FPP problems on the Configuration Model (CM), a model of constructing random networks with arbitrary degree distributions.
We shall defer a formal definition of this model to Section \ref{sec-not} and shall discuss
related work in Section \ref{sec-rlt}. Let it suffice to say that this model has arisen in myriad
applied contexts, ranging from combinatorics, computer science, statistical physics,
and epidemiology and seems to be one of the most widely used models in the modern
networking community.

We shall consider FPP on the CM where the  exponent $\tau$ of the degree distribution satisfies $\tau \in [1,2)$
and each edge is given a random exponential edge weight. FPP for
the case $\tau> 2$ was analyzed in \cite{BhaHofHoo09b} where the hopcount seems to exhibit a
remarkably universal behavior. More precisely, the hopcount always scales as $\log{n}$
and central limit theorems (CLTs) with matching asymptotic means and variances hold.
While these graphs are sparse and locally tree-like,  what is remarkable is that the same
fact also holds in the case of the most well-connected graph, namely the complete graph,
for which the hopcount satisfies a CLT as the model with $\tau>2$, with asymptotic mean and variance equal to
$\log{n}$, as $n\to\infty$. See, e.g., \cite{janson123} and \cite{vcg-random-shanky}
and the references therein.

When the degree exponent $\tau$ is in the interval $[1,2)$ we shall find that CLTs do {\it not} hold, and that
the hopcount remains uniformly bounded due to the remarkable shape of such networks, which we may
think of as a collection of interconnected star networks, the centers of the stars corresponding
to the vertices with highest degrees. We shall consider two models of the network topology,
one where we look at the original CM and the second more realistic model
called the \emph{erased} model where we shall delete all self loops and merge all multiple edges
from the original CM. In the resulting graph, each edge receives an independent
exponential weight with rate 1. Thus, for the erased CM, the direct weight between
two vertices connected by an edge is an exponential random variable with
rate 1, while for the original CM, it is also exponential, but with rate
equal to the number of edges between the pair of vertices. When $\tau>2$, there is no essential
difference between the original and the CM \cite{BhaHofHoo09b}.

In both cases, we shall see that the hopcount is tight and that a limit distribution exists.
More surprisingly, in the erased CM, this limiting distribution puts mass only on the {\it even} integers.
We also exhibit a nice constructive picture of how this arises, which uses the powerful machinery
of Poisson-Dirichlet distributions. We further find the distributional limit of the weight of
the minimal weight path joining two typical vertices.

Since the hopcount remains tight, this model is remarkably efficient in transporting
or routing flow between vertices in the network. However, a downside of
this property of the network is its extreme fragility w.r.t. directed attacks
on the network. More precisely, we shall show that there exists a simple algorithm deleting
a bounded number of
vertices such that the chance of disconnecting any two typical vertices is close to $1$ as
$n\to\infty$. At the same time we shall also show that these networks are relatively
stable against random attacks.


This paper is organized as follows. In Section \ref{sec-not}, we shall introduce the model and some
notation. In Section \ref{sec-results}, we state our main results. In Section \ref{sec-rlt},
we describe connections to the literature and discuss our results. In Section \ref{sec-proof(1,2)},
we give the proof in the original CM, and in Section \ref{sec:proofs-erased}, we prove the
results in the erased CM.

\section{Notation and definitions}
\label{sec-not}
In this section, we introduce the random graph model that we shall be working on, and recall some limiting
results on i.i.d.\ random variables with infinite mean. We shall use the notation that
$f(n)=O(g(n))$, as $n\rightarrow \infty$, if $|f(n)| \leq C g(n)$, and
$f(n)=o(g(n))$, as $n\rightarrow \infty$, if $|f(n)|/g(n)\rightarrow 0$.
For two sequences of random variables $X_n$ and $Y_n$,
we write that $X_n=\Op(Y_n)$, as $n\rightarrow \infty$,
when $\{X_n/Y_n\}_{n\geq 1}$ is a tight sequence of random variables.
We further write that $X_n=\Theta_{\sss \prob}(Y_n)$ if $X_n=\Op(Y_n)$ and $Y_n=\Op(X_n)$.
Further, we write that $X_n=\op(Y_n)$, when  $|X_n|/Y_n$ goes to $0$ in probability ($|X_n|/Y_n\convp 0$);
equality in distribution is denoted by the symbol $\sim$.
Throughout this paper, for a sequence of events $\{F_n\}_{n\geq 1}$,
we say say that $F_n$ \emph{occurs with high probability} ({\bf whp})
if $\lim_{n\rightarrow \infty} \prob(F_n)=1$.

\paragraph{Graphs:} We shall typically be working with random graphs on $n$ vertices, which
have a giant component consisting of $n-o(n)$ vertices. Edges are given a
random edge weight (sometimes alternatively referred to as cost)
which in this study will always be assumed to be independent,
exponentially distributed random variables with mean 1.
We pick two vertices uniformly at random in the network.
We let $W_n$ be the random variable denoting the total weight of the minimum weight path
between the two typical vertices and $H_n$ be the number of edges on this path or {\it hopcount}.

\paragraph{Construction of the configuration model:} We are interested in constructing a random graph on
$n$ vertices. Given  a {\bf degree sequence}, namely a sequence of $n$ positive integers
$\bfD = (D_1,D_2,\ldots, D_n)$ with the total degree
    \eqn{
    \lbeq{Ln-def}
    L_n=\sum_{i=1}^n D_i
    }
assumed to be even, the CM on $n$ vertices with
degree sequence $\bfD$ is constructed as follows:

Start with $n$ vertices and $D_j$ stubs adjacent to vertex $j$.
The graph is constructed by pairing up each stub to some other stub to form edges.
Number the stubs from $1$ to $L_n$ in some arbitrary order. Then, at each step,
two stubs (not already paired) are chosen uniformly at random among all the {\it free}
stubs and are paired to form a single edge in the graph.
These stubs are no longer free and removed from the list of free stubs.
We continue with this procedure of choosing and pairing two stubs until all the stubs are
paired.

\paragraph{Degree distribution:} The above denoted the construction of the CM when the
degree distribution is given and the total degree is even. Here we specify how we construct the actual degree sequence $\bfD$.
We shall assume that each of the random variables $D_1,D_2,\ldots D_n$ are independent
and identically distributed (i.i.d.) with distribution $F$.
(Note that if the sum of stubs $L_n$ is not even then we use the degree sequence with
$D_n$ replaced with $D_n+1$. This will not effect our calculations).

We shall assume that the degree distribution $F$, \rr{with atoms $f_1,f_2,\ldots$}
satisfies the property:
    \begin{equation}
    \label{defF}
     1-F(x) =x^{-(\tau-1)} L(x),
     \end{equation}
for some slowly varying function $x\mapsto L(x)$. Here, the parameter $\tau$,
which we shall refer to as the {\it degree exponent},
is assumed to be in the interval $[1,2)$, so that $\expec[D_i]=\infty$.
In some cases, we shall make stronger assumptions than \eqref{defF}.

\paragraph{Original model:} We assign to each edge a
random and i.i.d.\ exponential mean one edge weight.
Throughout the sequel, the weighted random graph so generated will be referred
to as the {\it original} model and we shall denote the random network so
obtained as $\GGor_n$.

\paragraph{Erased model:} This model is constructed as follows: Generate a CM as before
and then erase all self loops and merge all multiple edges into a single edge.
\emph{After} this, we put independent exponential weights with rate $1$ on the (remaining) edges.
Thus, while the graph distances are not affected by the erasure, we shall see that
the hopcount has a different limiting distribution. We shall denote the random
network on $n$ vertices so obtained by $\GGer_n$.

\subsection{Poisson-Dirichlet distribution}
\label{sec-PDdistr}
Before describing our results, we shall need to make a brief detour into extreme
value theory for heavy-tailed random variables. As in \cite{hofs2}, where the graph distances in
the CM with $\tau\in [1,2)$ are studied, the relative sizes of the order statistics of the
degrees play a crucial role in the proof. In order to describe the limiting behavior of
the order statistics, we need some definitions.

We define a (random) probability distribution $P=\{P_i\}_{i\geq 1}$ as
follows. Let $\{E_i\}_{i=1}^\infty$ be i.i.d.\ exponential random variables
with rate 1, and define $\Gamma_i=\sum_{j=1}^i E_j$.
Let $\{D_i\}_{i=1}^{\infty}$ be an i.i.d.\ sequence of
random variables with distribution function $F$ in \eqref{defF}, and let
$D_{\sss(n:n)}\geq D_{\sss(n-1:n)}\geq \cdots \geq D_{\sss(1:n)}$
be the order statistics of $\{D_i\}_{i=1}^{n}$. In the sequel of this paper,
we shall label vertices according to their degree, so that
vertex 1 has maximal degree, etc.

We recall \cite[Lemma 2.1]{hofs2}, that
there exists a sequence $u_n$, with $u_n= n^{1/(\tau-1)}l(n)$, where $l$ is
slowly varying, such that
    \eqan{
    \lbeq{EVTexpsrep}
    u_n^{-1} \left(L_n,\{D_{\sss(n+1-i:n)}\}_{i=1}^{\infty}\right)\convd
    \left(
    \sum_{j=1}^\infty \Gamma_j^{-1/(\tau-1)},
    \{\Gamma_i^{-1/(\tau-1)}\}_{i=1}^{\infty}
    \right),
    }
where $\convd$ denotes convergence in distribution.
We abbreviate $\xi_i=\Gamma_i^{-1/(\tau-1)}$ and
$\eta=\sum_{j=1}^\infty \xi_j$ and let
    \eqn{
    \lbeq{PRPDdef}
    P_i=\xi_i/\eta,\qquad i\geq 1,
    }
so that, $P=\{P_i\}_{i\geq 1}$
is a {\it random} probability distribution. The sequence $\{P_i\}_{i\geq 1}$ is called
the {\it Poisson-Dirichlet distribution} (see e.g., \cite{PitYor97}). A lot is known
about the probability distribution $P$. For example,
\cite[Eqn. (6)]{PitYor97} proves that for any $f\colon [0,1]\rightarrow {\mathbb R}$,
and with $\alpha=\tau-1\in (0,1)$,
    \eqn{
    \lbeq{mean-PDD}
    \expec\big[\sum_{i=1}^{\infty} f(P_i)\big]=
    \frac{1}{\Gamma(\alpha)\Gamma(1-\alpha)}\int_0^1 f(u) u^{-\alpha-1}(1-u)^{\alpha-1}du.
    }
For example, this implies that
    \eqn{
    \lbeq{p2-reason}
    \expec\big[\sum_{i=1}^{\infty} P_i^2\big]=\frac{\Gamma(\alpha)\Gamma(2-\alpha)}
    {\Gamma(\alpha)\Gamma(1-\alpha)}=1-\alpha=2-\tau.
    }

\section{Results}
\label{sec-results}
In this section, we state the main results of the paper, separating between the
original CM and the erased CM.

\subsection{Analysis of shortest-weight paths for the original CM}
\label{sec-tau(1,2)}

Before describing the results we shall need to construct a limiting infinite
object $\KKor$ in terms of the Poisson-Dirichlet distribution $\{P_i\}_{i\geq 1}$
given in \eqref{PRPDdef} and the sequence of random variables $\xi_i$ and their sum
$\eta$ which arise in the representation of this distribution.  This will be an
infinite graph with weighted edges on the vertex set $\Zbold^+ = \{1,2,\ldots\}$,
where every pair of vertices $(i,j)$ is connected by an edge which, conditionally
on $\{\xi_i\}_{i\geq 1}$, are independent exponential random variables with exponential
distribution with rate $\xi_i \xi_j/\eta$.

Let $\Wor_{ij}$ and $\Hor_{ij}$ denote the weight and number of edges
of the minimal-weight path in $\KKor$ between the vertices $i,j\in \Zbold^+$. Our
results will show that, in fact, the FPP problem on $\KKor$
is well defined \ch{(see Proposition \ref{prop-orig-well-defined}}.
Let $\Ior$ and $\Jor$ be two vertices chosen independently at
random from the vertex set $\Zbold^+$ with probability $\{P_i\}_{i\geq 1}$.
Finally, recall that $\GGor_n$ is the random network on $n$ vertices with
exponential edge weights constructed in Section \ref{sec-not}. We are now in a
position to describe our limiting results for the original CM:

\begin{Theorem}[Asymptotics FPP for the original CM]
\label{main(1,2)}
\ch{Consider the random network $\GGor_n$,  with the degree distribution
$F$ satisfying \eqref{defF} for some $\tau\in [1,2)$.}\\
(a) \ch{Let $\Wor_n$ be the weight of the minimal weight path between
two uniformly chosen vertices in the network. Then,}
    \eqn{
    \Wor_n \convd \Vor_1+\Vor_2,
    }
where $\Vor_i$, $i=1,2$, are independent random variables
with $\Vor_i\sim E_i/D_i$, where $E_i$ is exponential with rate 1 and
$D_1, D_2$ are independent and identically distributed with distribution $F$,
independently of $E_1, E_2$. More precisely, as $n\to\infty$,
    \eqn{
    u_n\big(\Wor_n-(\Vor_1+\Vor_2)\big)\convd \Wor_{\Ior \Jor},
    }
where $u_n$ is defined by
    \eqn{
    \lbeq{un-def}
    u_n=\sup\{u: 1-F(u)\geq 1/n\}.
    }
(b) Let $\Hor_n$ be the number of edges in the minimal weight path between
two uniformly chosen vertices in the network. Then,
    \eqn{
    \Hor_n \convd 2+ \Hor_{\Ior \Jor}.
    }
Writing $\pi_k = \prob(\Hor_{\Ior \Jor} = k-2)$, we have $\pi_k>0$ for each $k\geq 2$, when $\tau\in (1,2)$.
The probability distribution $\pi$ depends only on $\tau$, and not on any other
detail of the degree distribution $F$. Moreover,
    \eqn{
    \lbeq{pi2-form}
    \pi_2=2-\tau.
    }
\end{Theorem}

Theorem \ref{main(1,2)} implies that, for $\tau\in [1,2)$, the
hopcount is uniformly bounded, as is the case for the typical graph
distance obtained by taking the weights to be equal to 1 a.s.\ (see \cite{hofs2}).
However, while for unit edge weights and $\tau\in (1,2)$, the limiting hopcount is at most
3, for i.i.d.\ exponential weights the limiting hopcount can take all integer
values greater than or equal to 2.

\subsection{Analysis of shortest-weight paths for the erased CM}
\label{sec-erased}
The results in the erased CM hold under a more restricted condition on the degree
distribution $F$. More precisely, we assume that there exists a constant $0<c<\infty,$ such that
    \eqn{
    \lbeq{assumption-F-erased}
     1-F(x)= c x^{-(\tau-1)} (1+o(1)), \qquad x\rightarrow \infty,
    }
and we shall often make use of the upper bound $1-F(x)\leq c_2 x^{-(\tau-1)},$ valid for all $x\geq 0$
and some constant $c_2>0$.

Before we can describe our limit result for the erased CM,
we shall need an explicit construction of a limiting infinite network $\KKer$ using
the Poisson-Dirichlet distribution described in \eqref{PRPDdef}. Fix a realization
$\{P_i\}_{i\geq 1}$. Conditional on this sequence, let $f(P_i, P_j)$ be the probability
    \eqn{
    \label{eq:fpij}
    f(P_i, P_j) = \prob(\mathcal{E}_{ij}),
    }
of the following event $\mathcal{E}_{ij}$:
    \begin{quote}
    Generate a random variable $D\sim F$ where $F$ is the degree distribution. Conduct $D$
    independent multinomial trials where we select cell $i$ with probability $P_i$ at each stage.
    Then $\mathcal{E}_{ij}$ is the event that both cells $i$ and $j$ are \rr{selected}.
    \end{quote}
More precisely, for $0\le s,t\le 1$,
    \eqn{
    \lbeq{fpipj-def}
    f(s, t)=1-\expec[(1-s)^D]-\expec[(1-t)^D]+\expec[(1-s-t)^D].
    }

Now consider the following construction $\KKer$ of a random network on the vertex set $\Zbold^+$, where every vertex is connected to every other vertex by a single edge.  Further, each edge $(i,j)$ has a random weight $l_{ij}$
where, given $\{P_i\}_{i\geq 1}$, the collection $\{l_{ij}\}_{1\leq i<j<\infty}$ are conditionally
independent with distribution:
    \eqn{
    \label{eq:lij}
    \prob\left(l_{ij}>  x\right) = \exp\left(-f(P_i, P_j)x^2/2\right).
    }
Let $\Wer_{ij}$ and $\Her_{ij}$ denote the weight and number of edges
of the minimal-weight path in $\KKer$ between the vertices $i,j\in \Zbold^+$.
Our analysis shall, in particular, show that the FPP on $\KKer$ is well defined
(see Proposition \ref{prop-infi-FPP-Ker-well-defined} ).

Finally, construct the random variables $\Der$ and $\Ier$ as follows: Let $D\sim F$ and
consider a multinomial experiment with $D$ independent trials where at each trial,  we choose
cell $i$ with probability $P_i$. Let $\Der$ be the number of \emph{distinct}
cells so chosen and suppose the cells chosen are $\mathcal{A} = \{a_1, a_2, \ldots, a_{\Der}\}$.
Then let $\Ier$ be a cell chosen uniformly at random amongst $\mathcal{A}$.
Now we are in a position to describe the limiting distribution of the
hopcount in the erased CM:

\begin{Theorem}[Asymptotics FPP for the erased CM]
\label{thm-erased}
Consider the random network $\GGer_n$,  with the degree distribution
$F$ satisfying \refeq{assumption-F-erased} for some $\tau\in (1,2)$.\\
(a) Let $\Wer_n$ be the weight of the minimal weight path between
two uniformly chosen vertices in the network. Then,
    \eqn{
    \lbeq{Wn-weak-conv-simple}
    \Wer_n \convd \Ver_1 + \Ver_2.
    }
where $\Ver_i$, $i=1,2$, are independent random variables
with $\Ver_i\sim E_i/\Der_i$, where $E_i$ is exponential with rate 1 and
$\Der_1, \Der_2$ are, conditionally on $\{P_i\}_{i\geq 1}$,
independent random variables distributed as $\Der$, independently of $E_1, E_2$.
More precisely, as $n\to\infty$,
    \eqn{
    \lbeq{Wn-weak-conv-erased}
    \sqrt{n}\left(\Wer_n - (\Ver_1 + \Ver_2) \right) \convd \Wer_{\Ier\Jer}.
    }
(b) Let $\Her_n$ be the number of edges in the minimal weight path between
two uniformly chosen vertices in the network. Then,
    \eqn{
    \lbeq{Hn-weak-conv-erased}
    \Her_n \convd 2+ 2\Her_{\Ier \Jer},
    }
where $\Ier,\Jer$ are two copies of the random variable $\Ier$ described above,
which are conditionally independent given $P=\{P_i\}_{i\geq 1}$.
In particular, the limiting probability measure of the hopcount is
supported only on the even integers.
\end{Theorem}

We shall now present an intuitive explanation of the results claimed in Theorem \ref{thm-erased},
starting with \refeq{Wn-weak-conv-simple}. We let $\UnVer_1$ and $\UnVer_2$ denote
two uniformly
chosen vertices, note that they can be identical with probability $1/n$.
We further note that both vertex $\UnVer_1$ and $\UnVer_2$ have a random degree
which are close to independent copies of $D$.
We shall informally refer to the vertices with degrees $\Theta_{\sss \prob}(n^{1/(\tau-1)})$
as \emph{super vertices} (see \refeq{super-vertices-def} for a precise definition, and
recall \refeq{EVTexpsrep}).
We shall frequently make use of the fact that normal vertices are, {\bf whp},
exclusively attached to super vertices.
The number of super vertices to which $\UnVer_i,\, i=1,2,$ is attached to is equal
to $\Der_i,\,i=1,2$, as described above. The minimal weight edge between
$\UnVer_i,\, i=1,2,$ and any of its neighbors is hence equal in distribution
to the minimum of a total of $\Der_i$ independent exponentially distributed
random variables with mean 1. The shortest-weight path between two super vertices can pass through
intermediate normal vertices, of which there are $\Theta_{\sss \prob}(n)$.
This induces that the minimal weight between any pair of super vertices is of
order $\op(1)$, so that the main contribution to $\Wer_n$ in \refeq{Wn-weak-conv-simple}
is from the two minimal edges coming out of the vertices $\UnVer_i, i=1,2$.
This shows \refeq{Wn-weak-conv-simple} on an intuitive level.

We proceed with the intuitive explanation of \refeq{Wn-weak-conv-erased}. We use that,
{\bf whp}, the vertices $\UnVer_i,\, i=1,2,$ are only attached to super
vertices. Thus, in \refeq{Wn-weak-conv-erased},
we investigate the shortest-weight paths between super vertices.
Observe that we deal with the erased CM, so between any pair of
vertices there exists only {\it one} edge having an exponentially distributed
weight with mean 1. As before, we number the super vertices
by $i=1,2,\ldots$ starting from the largest degree.
We denote by $\Ner_{ij}$, the number of common neighbors of
the super vertices $i$ and $j$, for which we shall show that
$\Ner_{ij}$ is $\Theta_{\sss \prob}(n)$.

Each element in $\Ner_{ij}$ corresponds to a unique two-edge path
between the super vertices $i$ and $j$. Therefore, the weight of the
minimal two-edge path between the super vertices $i$ and $j$ has distribution
$w_{ij}^{\sss(n)}\equiv \min_{s\in \Ner_{ij}}(E_{is}+E_{sj}).$
Note that $\{E_{is}+E_{sj}\}_{s\in \Ner_{ij}}$ is a collection of
$\Ner_{ij}$ i.i.d.\ Gamma(2,1) random variables.
More precisely, $\Ner_{ij}$ behaves as $nf(P_i^{\sss(n)},P_j^{\sss(n)})$, where
$P_i^{\sss(n)}=D_{\sss (n+1-i:n)}/L_n$. Indeed, when we consider an \emph{arbitrary} vertex with degree
$D\sim F$, the conditional probability, conditionally on $\{P_i^{\sss(n)}\}_{i\geq 1}$,
that this vertex is both connected to super vertex $i$ and super vertex $j$ equals
    $$
    1-(1-P_i^{\sss(n)})^D-(1-P_j^{\sss(n)})^D
    +(1-P_i^{\sss(n)}-P_j^{\sss(n)})^D.
    $$
Thus, the expected number of vertices connected to both super vertices $i$ and $j$
is, conditionally on $\{P_i^{\sss(n)}\}_{i\geq 1}$,
$\Ner_{ij}\approx nf(P_i^{\sss(n)},P_j^{\sss(n)}),$
and $f(P_i^{\sss(n)},P_j^{\sss(n)})$ weakly converges to $f(P_i,P_j)$.

We conclude that the minimal two-edge path between super vertex $i$ and
super vertex $j$ is the minimum of $n f(P_i^{\sss(n)},P_j^{\sss(n)})$
Gamma(2,1) random variables $Y_s$, which are close to being independent.
Since
    \eqn{
    \lbeq{extrem-value}
    \lim_{n\rightarrow \infty} \prob(\sqrt{n}\min_{1\le s \le \beta n} Y_s>x)={\mathrm e}^{-\beta x^2/2},
    }
for any $\beta>0$,  \refeq{extrem-value}, with
$\beta=\beta_{ij}=f(P_i^{\sss(n)},P_j^{\sss(n)})\approx f(P_i,P_j)$ explains
the weights $l_{ij}$ defined in \refeq{eq:lij}, and also explains
intuitively why \refeq{Wn-weak-conv-erased} holds.

The convergence in \refeq{Hn-weak-conv-erased} is explained in a
similar way. Observe that in \refeq{Hn-weak-conv-erased}
the first $2$ on the right side originates from the 2 edges that connect $\UnVer_1$
and $\UnVer_2$ to the minimal-weight super vertex. Further, the factor $2$ in front
of $\Her_n$ is due to the fact that shortest-weight paths between super
vertices are concatenations of two-edge paths with random weights
$l_{ij}$. We shall further show that two-edge paths, consisting of an alternate
sequence of super and normal vertices, are the optimal paths in the sense of
minimal weight paths between super vertices.

This completes the intuitive explanation of Theorem \ref{thm-erased}.

\subsection{Robustness and fragility}
\label{sec-rob-frag}
The above results show that the hopcount $H_n$ in both models converges
in distribution as $n\to\infty$. Interpreting the hopcount
as the amount of travel time it takes for messages to get
from one typical vertex to another typical vertex, the above shows that the CM with $\tau\in (1,2)$
is remarkably efficient in routing flow between vertices.
We shall now show that there exists a down side to this
efficiency. The theorem is stated for the more natural erased CM but
one could formulate a corresponding theorem for the original CM as well.

\begin{Theorem}[Robustness and fragility]
\label{thm-rob-frag}
Consider the random weighted network $\GGer_n,$ where the degree distribution satisfies
\refeq{assumption-F-erased} for some $\tau\in (1,2).$ Then, the following properties hold:
\\(a) {\bf Robustness:} Suppose an adversary attacks the network via
randomly and independently deleting each vertex with probability $1-p$ and leaving each vertex with
probability $p$. Then, for any $p>0$, there exists a unique giant
component of size $\Theta_{\sss \prob}(n)$.\\
(b) {\bf Fragility:} Suppose an adversary attacks the network via deleting vertices of
maximal degree. Then, for any $\eps>0$, there exists an integer $K_\eps<\infty$ such that deleting
the $K_\eps$ maximal degree vertices implies that, for two vertices
$A_1$ and $A_2$ chosen uniformly at random from $\GGer_n$,
    \eqn{
    \lbeq{disconnection-prob}
    \limsup_{n\to\infty} \prob\left(A_1 \leftrightarrow A_2\right) \leq \eps.
    }
where $A_1\leftrightarrow A_2$ means that there exists a path connecting vertex $A_1$ and $A_2$ after
deletion of the maximal vertices. Thus one can disconnect the network by deleting \ch{$\Op(1)$}
vertices.
\end{Theorem}

\paragraph{Remark:} As in much of percolation theory, one could ask for the size of the
giant component in  part (a) above when we randomly delete vertices. See Section
\ref{sec-rob-frag-pf},  where we find the size of the giant component as $n\to\infty$,
and give the idea of the proofs for the reported behavior.


\section{Discussion and related literature}
\label{sec-rlt}
\ch{In this section, we discuss the literature and state some
further open problems and conjectures.}

\paragraph{The configuration model.}
The CM was introduced by Bender and Canfield \cite{BenCan78}, see also
Bollob\'as \cite{Boll01}. Molloy and Reed \cite{MolRee95} were the first
to use specified degree sequences. The model has become quite popular and
has been used in a number of diverse fields. See in particular
\cite{babak-newman2, babak-newman1} for applications to modeling of
disease epidemics and \cite{newman} for a full survey of various questions
from statistical physics.

For the CM, the graph distance, i.e., the minimal number of edges on a path
connecting two given vertices, is well understood. We refer to \cite{hofs3} for $\tau>3$, \cite{hofs1,NorRei04}
for $\tau\in (2,3)$ and \cite{hofs2} for $\tau \in (1,2)$. In the latter paper, it was shown that
the graph distance weakly converges, where the limit is either two or three, each with positive
probability.

\paragraph{FPP on random graphs.}
Analysis of FPP in the context of modern random graph models
has started only recently (see \cite{vcg-random-shanky,hofs-erdos-fpp,hofs-flood,janson123,wastlund}).
The particular case of the CM with degree distribution $1-F(x)=x^{1-\tau}L(x)$, where $\tau>2$, was studied in  \cite{BhaHofHoo09b}.
For $\tau>2$, where , the hopcount remarkably scales as
$\Theta(\log{n})$ and satisfies a central limit theorem (CLT) with asymptotic mean and
variance both equal to $\alpha\log{n}$ for some $\alpha>0$ (see \cite{BhaHofHoo09b}), this
despite the fact that for $\tau\in (2,3)$, the graph distance scales as $\log{\log{n}}$.
The parameter $\alpha$ belongs to $(0,1)$ for $\tau\in (2,3)$, while $\alpha>1$ for
$\tau>3$ and is the only feature which is left over from
the randomness of the random graph. As stated in Theorem \ref{main(1,2)} and \ref{thm-erased}, the behavior for $\tau\in (1,2)$,
where the hopcount remains bounded and weakly converges, is rather different from
the one for $\tau>2$.


\paragraph{Universality of $\KKor$ and $\KKer$.} Although we have used exponential edge weights,
we believe that one obtains the same result with any ``similar'' edge weight distribution
with a density $g$ satisfying $g(0) = 1$. More precisely, the hopcount result,
the description of $\KKor$ and $\KKer$ and the corresponding limiting distributions
in Theorems \ref{main(1,2)}--\ref{thm-erased} will remain unchanged. The only
thing that will change is the distribution of $(\Vor_1, \Vor_2)$ and $(\Ver_1, \Ver_2)$.
In Section \ref{sec-conc}, Theorem 8.1, we state what happens when the weight density $g$ satisfies $g(0)=\zeta\in (0,\infty)$.
When the edge weight density
$g$ satisfies $g(0) = 0$ or $g(0)=\infty$, then we expect that the hopcount
remains tight, but that the weight of the minimal path $W_n$, as well as the
limiting FPP problems, both for the original and erased CM, are different.

\paragraph{Robustness and fragility of random networks.}
The issue of robustness, yet fragility, of random network models has stimulated
an enormous amount of research in the recent past. See \cite{ba-rob} for one of the
original statistical physics papers on this topic, and \cite{boll-rior-perc}
for a rigorous derivation of this fact when the power-law exponent $\tau = 3$ in the
case of the preferential attachment model. The following universal property is
believed to hold for a wide range of models:
\begin{quote}
If the degree exponent $\tau$ of the model is in $(1,3]$, then the network is
robust against random attacks but fragile against directed attacks, while for
$\tau> 3$, under random deletion of vertices  there exists a critical
(model dependent $p_c$) such that for $p< p_c$ there is no giant component,
while for $p> p_c$, there is a giant component.
\end{quote}
Proving these results in a wide degree of generality is a challenging program in modern applied probability.

\paragraph{Load distributions on random networks.}
Understanding the FPP model on these networks opens
the door to the analysis of more complicated functionals such as the load distribution
on various vertices and edges  of the network, which measure the ability of the
network in dealing with congestion when transporting material from one part of the
network to another. We shall discuss such questions \ch{in some more detail}
in Section \ref{sec-conc}.

\paragraph{Organization of the proofs and conventions on notation.}
\label{sec-proofs}
The proofs in this paper are organized as follows.
In Section \ref{sec-proof(1,2)} we prove the
results for the original CM, while Section \ref{sec:proofs-erased} contains
the proofs for the erased CM. Theorem \ref{thm-rob-frag} is proved in
Section \ref{sec-rob-frag-pf}, and we close with a conclusion and discussion
in Section \ref{sec-conc}.

In order to simplify notation, we shall drop
the superscripts ${\tt er}$ and ${\tt or}$ so that for example the minimal
weight random variable $\Wor_n$ between two uniformly selected vertices will
be denoted by $W_n$ when proving facts about the original CM in
Section \ref{sec-proof(1,2)}, while $W_n$ will be used to denote $\Wer_n$
when proving facts about the erased CM in Section \ref{sec:proofs-erased}.

\section{Proofs in the original CM: Theorem \ref{main(1,2)} }
\label{sec-proof(1,2)}
In this section, we prove Theorem \ref{main(1,2)}. As part of the proof, we also prove
that the FPP on $\KKor$ is well defined, as formalized in the following
proposition:

\begin{Proposition}[FPP on $\KKor$ is well defined]
\label{prop-orig-well-defined}
For any fixed $K\geq 1$ and for all $i,j< K$ in $\KKor$, we have $\Wor_{ij}> 0$ for $i\neq j$
and $\Hor_{ij}< \infty$. In particular, this implies that $\Hor_{\Ior\Jor}<\infty$
almost surely, where we recall that $\Ior$ and $\Jor$
are two random vertices in $\Zbold_+$ chosen (conditionally) independently with
distribution $\{P_i\}_{i\geq 1}$.
\end{Proposition}

Recall that we label vertices according to their degree.
We let $\UnVer_1$ and $\UnVer_2$ denote two uniformly chosen vertices.
\ch{Since the CM has a giant component containing $n-o(n)$ vertices,
{\bf whp}, $\UnVer_1$ and $\UnVer_2$ will be connected.}
We note that the edge incident to vertex \rr{$\UnVer_1$}
with minimal weight has weight given by $V_i=E_i/D_{\UnVer_i},\, i=1,2$,
where $D_{\UnVer_1}$ denotes the degree of vertex $\UnVer_1$.
As a result, $(V_1, V_2)$ has the same distribution as
$(E_1/D_1, E_2/D_2)$, where $(D_1,D_2)$ are two independent random
variables with distribution function $F$.
Further, by \cite[Theorem 1.1]{hofs2}, \rr{{\bf whp}}, the
vertices $\UnVer_1$ and $\UnVer_2$ are not directly
connected.  When $\UnVer_1$ and $\UnVer_2$ are not directly
connected, then $W_n \geq V_1+V_2$, and $V_1$ and $V_2$ are
independent, as they depend on the exponential weights of {\it
disjoint} sets of edges, while, by construction, $D_{\UnVer_1}$ and
$D_{\UnVer_2}$ are independent. This proves the required lower bound in
Theorem \ref{main(1,2)}(a). For the upper bound, we further note
that, by \cite[Lemma 2.2]{hofs2}, the vertices $\UnVer_1$ and $\UnVer_2$
are, {\bf whp}, exclusively connected to so-called
{\it super vertices}, which are the $m_n$ vertices with the
largest degrees, for any $m_n\rightarrow \infty$ arbitrarily
slowly. Thus, the upper bound follows if any two of such super
vertices are connected by an edge with weight which converges to 0
in distribution. Denote by $M_{i,j}$ the minimal weight of all
edges connecting the vertices $i$ and $j$. Then, conditionally on
the number of edges between $i$ and $j$, we have that $M_{i,j}\sim
{\rm Exp}(N(i,j))$, where $N(i,j)$ denotes the number of edges
between $i$ and $j$, and where we use ${\rm Exp}(\lambda)$ to denote an exponential random variable with \rr{rate $\lambda$}.
We further denote
$P_i^{\sss(n)}=D_{\sss(n+1-i:n)}/L_n$, so that
$P^{\sss(n)}=\{P_i^{\sss(n)}\}_{i=1}^n$ converges in distribution
to the Poisson-Dirichlet distribution. \rr{We will} show
that, conditionally on the degrees and {\bf whp},
    \eqn{
    \lbeq{Nij-asy}
    N(i,j)=(1+\op(1)) L_n P_i^{\sss(n)}P_j^{\sss(n)}.
    }
Indeed, we note that
    \eqn{
    N(i,j)=\sum_{s=1}^{D_i} I_{s}(i,j),
    }
where $I_{s}(i,j)$ is the indicator that the $s^{\rm th}$ stub
of vertex $i$ connects to $j$. We write $\prob_n$ for the conditional
distribution given the degrees, and $\expec_n$ for the expectation
w.r.t.\ $\prob_n$. It turns out that we can even
prove Theorem \ref{main(1,2)} {\it conditionally on the degrees},
which is stronger than Theorem \ref{main(1,2)} {\it averaged}
over the degrees. For this, we note that,
for $1\leq s_1<s_2\leq D_i$,
    \eqn{
    \prob_n(I_{s_1}(i,j)=1)=\frac{D_j}{L_n-1},
    \qquad \prob_n(I_{s_1}(i,j)=I_{s_2}(i,j)=1)=\frac{D_j(D_j-1)}{(L_n-1)(L_n-3)},
    }
which implies, further using that $D_j=D_{\sss(n+1-j:n)}$ and
thus $D_j/L_n\convd P_j$, that
    \eqn{
    {\rm Var}_n(N(i,j))\leq C\frac{D_i^2 D_j}{L_n^2}=\op\Big(\frac{D_i^2D_j^2}{L_n^2}\Big)
    =\op\big(\expec_n[N(i,j)]^2\big).
    }
As a result, $N(i,j)$ is concentrated, and thus \refeq{Nij-asy}
follows.

In particular, we see that the vector
$\{N(i,j)/L_n\}_{i,j=1}^n$ converges in distribution
to $\{P_iP_j\}_{i,j=1}^{\infty}$. Thus,
for every $i,j$, and conditionally on the degrees,
we have that $M_{i,j}$ is approximately
equal to an exponential random variable with asymptotic
mean $L_n P_iP_j$. This proves that, with $J_1$ and $J_2$ being two random
variables, which are independent, conditionally on $P=\{P_i\}_{i=1}^\infty$, and with
    \eqn{
    \lbeq{Js-distr}
    \prob(J_s=i|P)=P_i,
    }
we have that
    \eqn{
    V_1+V_2\leq W_n \leq V_1+V_2+{\rm Exp}(L_n P_{J_1}P_{J_2}).
    }
Consequently, $u_n\big(W_n
-(V_1+V_2)\big)$ is a tight random variable. Below, we shall prove that, in fact,
$u_n\big(W_n-(V_1+V_2)\big)$ converges weakly to a non-trivial random variable.

Recall the above analysis, and recall that the edges with minimal
weight from the vertices $\UnVer_1$ and $\UnVer_2$ are connected to vertices
$J_1$ and $J_2$ with asymptotic probability, conditionally on the
degrees, given by \refeq{Js-distr}. Then, $H_n=2$ {\it precisely}
when $J_1=J_2$, which occurs, by the conditional independence of
$J_1$ and $J_2$ given $P$, with asymptotic probability
    \eqn{
    \prob_n(H_n=2)=\sum_{i=1}^{\infty} (P_i^{\sss(n)})^2 +\op(1).
    }
Taking expectations and using \refeq{p2-reason} together with the
bounded convergence theorem proves \refeq{pi2-form}.

Recall that $J_1$ and $J_2$ are
the vertices to which the edges with minimal weight
from $\UnVer_1$ and $\UnVer_2$ are connected, and recall their distribution in
\refeq{Js-distr}. We now prove the weak convergence of $H_n$ and of
$u_n\big(W_n-(V_1+V_2)\big)$ by constructing a shortest-weight tree in
$\KKor$.

We start building the shortest-weight tree from
$J_1$, terminating when $J_2$ appears for the first time in this tree.
We denote the tree of size $l$ by $T_l$, and note that
$T_1=\ch{\{J_1\}}$. Now we have the following recursive procedure to
describe the asymptotic distribution of $T_l$.
We note that, for any set of vertices $A$, the edge
with minimal weight outside of $A$ is a uniform edge
pointing outside of $A$.
When we have already constructed $T_{l-1}$, and we fix $i\in T_{l-1}, j\not\in T_{l-1}$,
then by \refeq{Nij-asy} there are approximately $L_n P_iP_j$
edges linking $i$ and $j$. \rr{Thus, the probability that vertex $j$ is added to $T_{l-1}$
is, conditionally on $P$, approximately equal to}
    \eqn{
    p_{ij}(l)
    =
    \frac{L_n P_j \sum_{a\in T_{l-1}} P_a}
    {L_n\sum_{a\in T_{l-1}, b\not\in T_{l-1}}  P_aP_b}
    =\frac{P_j}{1-P_{T_{l-1}}}\geq P_j,
    }
    where, for a set of vertices $A$, we write
    \eqn{
    P_A=\sum_{a\in A} P_a.
    }
Denote by $B_l$ the $l^{\rm th}$ vertex chosen.
We stop this procedure when $B_l=J_2$ for the first time,
and denote this stopping time by $S$,
so that, {\bf whp}, $H_n=2+H(S)$, where $H(S)$ is the height
of $B_S$ in  $T_{S}$. Also, $u_n\big(W_n-(V_1+V_2)\big)$ is
equal to $W_{S}$, which is the weight of
the path linking $J_1$ and $J_2$ in $\KKor$.

Note that the above procedure terminates in finite time, since $P_{J_2}>0$ and
at each time, we pick $J_2$ with probability at least $P_{J_2}$.
This proves that $H_n$ weakly converges,
and that the distribution is given only in terms of $P$.
Also, it proves that the FPP problem on $\KKor$ is well defined,
\ch{as formalized in Proposition \ref{prop-orig-well-defined}.}

\ch{Further,} since the distribution of $P$ only depends on $\tau\in [1,2)$,
and not on any other details of the degree distribution $F$,
the same follows for $H_n$. \ch{When $\tau=1$, then $P_1=1$ a.s.,
so that $\prob_n(H_n=2)=1+\op(1)$. When $\tau\in (1,2)$, on the other hand,
$P_i>0$ a.s.\ for each $i\in \N$, so that, by the above construction,
it is not hard to see that $\lim_{n\rightarrow \infty}
\prob_n(H_n=k)=\pi_k(P)>0$ a.s.\ for each $k\geq 2$.}
Thus, the same follows for $\pi_k=\lim_{n\rightarrow \infty} \prob(H_n=k)
=\expec[\pi_k(P)]$. It would be of interest to compute $\pi_k$ for $k>2$
explicitly, or even $\pi_3$, but this seems a difficult problem.
\qed

\section{Proofs in the erased CM: Theorem \ref{thm-erased}}
\label{sec:proofs-erased}
In this section, we prove the various results in the erased setup.
We start by giving an overview of the proof.

\subsection{Overview of the proof of Theorem \ref{thm-erased}}
\label{sec-overview-erased}
In this section, we formulate four key propositions, which, together, shall make the
intuitive proof given below Theorem \ref{thm-erased} precise, and which shall
combine to a formal proof of Theorem \ref{thm-erased}.

As before, we label vertices by their (original) degree so that vertex $i$ will be the
vertex with the $i^{\rm th}$ largest degree. Fix a sequence $\eps_n\to 0$
arbitrarily slowly. Then, we define the set of {\it super vertices}
$\su_n$ be the set of vertices with largest degrees, namely,
    \eqn{
    \lbeq{super-vertices-def}
    \su_n = \{i: D_i > \eps_n n^{1/(\tau-1)}\}.
    }
We shall refer to $\su_n^c$ as the set of {\it normal vertices}.

Recall the definition of the limiting infinite ``complete graph'' $\KKer$ defined in
Section \ref{sec-erased} and for any fixed $ k\geq 1$, let $(\KKer)^k$ denote the projection
of this object onto the first $k$ vertices (so that we retain only the first $k$ vertices
$1,2,\ldots, k$ and the corresponding edges between these vertices). Then the following
proposition says that we can move between the super vertices via two-edge paths which have
weight $\Theta(1/\sqrt{n})$. For notational convenience, we write $[k]:= \{1,2 \ldots, k\}$.

\begin{Proposition}[Weak convergence of FPP problem]
\label{prop-weak-con-FPP-erased}
Fix $k$ and consider the subgraph \rr{of the CM} formed by retaining the maximal $k$
vertices and all paths connecting any pair of these vertices by a single
intermediary normal vertex (i.e., two-edge paths). For any pair of vertices $i,j\in [k]$,
let $l_{ij}^{\sss (n)} = \sqrt{n} w_{ij}^{\sss (2)},$ where $w_{ij}^{\sss (2)}$ is the minimal
weight of all two-edge paths between $i$ and $j$ (with $w_{ij}^{\sss (2)} = \infty$
if they are not connected by a two-edge path).
Consider the complete graph $\KK_n^k$ on vertex set
$[k]$ with edge weights $l_{ij}^{\sss (n)}$. Then,
    \eqn{
    \lbeq{KK-conv}
    \KK_n^k \convd (\KKer)^k,
    }
where $\convd$ denotes the usual finite-dimensional
convergence of the ${k \choose 2}$ random variables $l_{ij}^{\sss (n)}$.
\end{Proposition}

The proof of Proposition \ref{prop-weak-con-FPP-erased} is deferred to Section
\ref{sec-pf-prop6.1}. Proposition \ref{prop-weak-con-FPP-erased} implies that the FPP problem
on the first $k$ super vertices along the two-edge paths converges in distribution
to the one on $\KKer$ restricted to $[k]$. We next
investigate the structure of the minimal weights from a uniform
vertex, and the tightness of recentered minimal weight:

\begin{Proposition}[Coupling of the minimal edges from uniform vertices]
\label{prop-coupling+tightness-weight}
Let $(\UnVer_1,\UnVer_2)$ be two uniform vertices, and let
$(V^{\sss(n)}_1, V^{\sss(n)}_2)$ denote the minimal weight in the erased CM
along the edges attached to $(\UnVer_1,\UnVer_2)$.\\
(a) Let $I^{\sss(n)}$ and $J^{\sss(n)}$ denote the vertices to which $\UnVer_i,\, i=1,2,$ are connected,
and let $(I,J)$ be two random variables having the distribution specified
right before Theorem \ref{thm-erased}, which are conditionally independent
given $\{P_i\}_{i\geq 1}$. Then, we can couple $(I^{\sss(n)},J^{\sss(n)})$
and $(I,J)$ in such a way that
    \eqn{
    \lbeq{coupling-IJ}
    \prob\big((I^{\sss(n)},J^{\sss(n)})\neq (I,J)\big)=o(1).
   }
(b) Let $V_i=E_i/\Der_i$, where $(\Der_1,\Der_2)$ are
two copies of the random variable $\Der$ described right before Theorem \ref{thm-erased},
which are conditionally independent given $\{P_i\}_{i\geq 1}$.\\
Then, we can couple $(V^{\sss(n)}_1, V^{\sss(n)}_2)$
to $(V_1,V_2)$ in such a way that
    \eqn{
    \lbeq{coupling-Vs}
    \prob\big((V^{\sss(n)}_1, V^{\sss(n)}_2)\neq (V_1,V_2)\big)=o(1).
    }
As a result, the recentered random variables $\sqrt{n}\big(W_n-(V_1+V_2)\big)$
form a tight sequence.
\end{Proposition}

The proof of Proposition \ref{prop-coupling+tightness-weight} is deferred to Section
\ref{sec-pf-prop6.2}. The following proposition asserts that
the hopcount and the recentered weight between the first $k$ super vertices
are tight random variables, and, in particular,
they remain within the first $[K]$ vertices, {\bf whp},
as $K\rightarrow \infty$:

\begin{Proposition}[Tightness of FPP problem and evenness of hopcount]
\label{prop-hopcount-conv}
Fix $k\geq 1$. For any pair of vertices $i,j\in [k]$, let $H_n(i,j)$
denote the number of edges of the minimal-weight path
between $i$ and $j$. Then,\\
(a) $H_n(i,j)$ is a tight sequence of random variables,
which is such that $\prob(H_n(i,j)\not\in 2\Zbold^+)=o(1)$;\\
(b) the probability that any of the minimal weight paths between $i,j\in [k]$, at even
times, leaves the $K$ vertices of largest degree tends to zero
when $K\rightarrow \infty$;\\
(c) the hopcount $H_n$ is a tight sequence of random variables,
which is such that $\prob(H_n\not\in 2\Zbold^+)=o(1)$.
\end{Proposition}

The proof of Proposition \ref{prop-hopcount-conv} is deferred to Section
\ref{sec-pf-prop6.3}. The statement is consistent with the intuitive explanation
given right after Theorem \ref{thm-erased}: the minimal weight paths between
two uniform vertices consists of an \emph{alternating} sequence of normal
vertices and super vertices. We finally state that the infinite FPP on the erased
CM is well defined:

\begin{Proposition}[Infinite FPP is well defined]
\label{prop-infi-FPP-Ker-well-defined}
Consider FPP on $\KKer$ with weights $\{l_{ij}\}_{1\leq i<j<\infty}$ defined in
\eqref{eq:lij}. Fix $k\geq 1$ and $i,j\in [k]$. Let ${\cal A}_{\sss K}$ be
the event that there exists a path
of weight at most $W$ connecting $i$ and $j$, which contains a
vertex in $\Zbold^+\setminus [K]$, and which is of weight at most $W$. Then,
there exists a $C>0$ such that, for all $K$ sufficiently large,
    \eqn{
    \lbeq{infi-FPP-bound}
    \prob({\cal A}_{\sss K})\leq CW K^{-1} {\mathrm e}^{CW\sqrt{\log{K}}}.
    }
\end{Proposition}

The proof of Proposition \ref{prop-infi-FPP-Ker-well-defined} is deferred to Section
\ref{sec-pf-prop6.4}. With Propositions \ref{prop-weak-con-FPP-erased}--\ref{prop-infi-FPP-Ker-well-defined}
at hand, we are able to prove Theorem \ref{thm-erased}:
\medskip

\noindent
{\it Proof of Theorem \ref{thm-erased} subject to Propositions \ref{prop-weak-con-FPP-erased}--\ref{prop-infi-FPP-Ker-well-defined}.}
By Proposition \ref{prop-coupling+tightness-weight}(b), we can couple
$(V^{\sss(n)}_1, V^{\sss(n)}_2)$ to $(V_1,V_2)$ in such a way that $(V^{\sss(n)}_1, V^{\sss(n)}_2)= (V_1,V_2)$ occurs {\bf whp}.
Further, {\bf whp}, for $k$ large, $I,J\leq k$, which
we shall assume from now on, while, by Proposition \ref{prop-coupling+tightness-weight}(b),
$\sqrt{n}\big(W_n-(V_1+V_2)\big)$ is a tight sequence
of random variables.

By Proposition \ref{prop-hopcount-conv}, the hopcount is a tight
sequences of random variables, which is {\bf whp} even. Indeed,
it consist of an alternating sequence of normal and super vertices.
We shall call the path of super vertices the \emph{two-edge path}. Then,
Proposition \ref{prop-hopcount-conv} implies that the probability that any
of the two-edge paths between any of the first $[k]$ vertices leaves the
first $K$ vertices is small when $K$ grows big.
As a result, we can write
$H_n=2+2H^{\sss(n)}_{I^{\sss(n)}J^{\sss(n)}}$, where $H^{\sss(n)}_{I^{\sss(n)}J^{\sss(n)}}$ is
the number of two-edge paths in $\KKner$. By \refeq{coupling-IJ}, we have that,
{\bf whp}, $H^{\sss(n)}_{I^{\sss(n)}J^{\sss(n)}}=H^{\sss(n)}_{IJ}$.

By Proposition \ref{prop-weak-con-FPP-erased},
the FPP on the $k$ vertices of largest degree in the CM weakly converges to the FPP
on the first $k$ vertices of $\KKer$, for any $k\geq 1$.
By Proposition \ref{prop-infi-FPP-Ker-well-defined}, {\bf whp}, the shortest-weight path
between any two vertices in $[k]$ in $\KKer$ does not leave the first $K$ vertices,
so that $W_{IJ}$ and $H_{IJ}$ are finite random variables, where
$W_{IJ}$ and $H_{IJ}$ denote the weight and number of steps in the
minimal path between $I$ and $J$ in $\KKer$. In particular, it follows
that $\sqrt{n}\big(W_n-(V^{\sss(n)}_1+V^{\sss(n)}_2)\big)\convd W_{IJ}$,
and that $H^{\sss(n)}_{ij}\convd H_{ij}$ for every $i,j\in [k]$, which is
the number of hops between $i,j\in [k]$ in $\KKer$. Since, {\bf whp},
$(V_1,V_2)=(V^{\sss(n)}_1,V^{\sss(n)}_2)$,
$\sqrt{n}\big(W_n-(V_1+V_2)\big)$ converges to the same limit.
This completes the proof
of Theorem \ref{thm-erased} subject to Propositions \ref{prop-weak-con-FPP-erased}--\ref{prop-infi-FPP-Ker-well-defined}.
\qed

\subsection{Weak convergence of the finite FPP problem to $\KKer$: Proof of Proposition \ref{prop-weak-con-FPP-erased}}
\label{sec-pf-prop6.1}
In this section, we study the weak convergence of the FPP on $\KK_n^k$ to the one on $(\KKer)^k$, by proving Proposition \ref{prop-weak-con-FPP-erased}.

We \rr{start by proving} some elementary results regarding the extrema of
Gamma random variables. We start with a particularly simple case, and
after this, generalize it to the convergence of all weights of two-edge paths in
$\KKner$.

\begin{Lemma}[Minima of Gamma random variables]
\label{lem-min-Gamma}
(a) Fix $\beta>0$ and  consider $n \beta$ i.i.d.\ Gamma(2,1) random variables $Y_i$.
Let $T_n=\min_{1\leq i\leq \beta n} Y_i$ be the
minimum of these random variables. Then, as $n\rightarrow \infty$,
    \eqn{
    \prob(\sqrt{n}T_n > x )
    \to \exp\left(-\beta x^2/2\right).
    \label{eq:lim-lij}}
(b) Let $\{X_i\}_{1\leq i\leq m}, \{Y_i\}_{1\leq i\leq m}\mbox{ and }\{Z_i\}_{1\leq i\leq m} $ be all independent collections of
independent exponential mean $1$ random variables. Let
    \eqn{
    \eta_m=\sqrt{m}\min_{1\leq i\leq m} (X_i+ Y_i),\qquad \kappa_m = \sqrt{m}\min_{1\leq i\leq m} (X_i+ Z_i),
    \qquad \text{and}\qquad \rho_m = \sqrt{m}\min_{1\leq i\leq m} (Y_i+ Z_i).
    }
Then, as $m\rightarrow \infty$,
    \eqn{
    (\eta_m, \kappa_m, \rho_m) \convd (\zeta_1, \zeta_2,\zeta_3).
    }
Here $\zeta_i$ are independent with the distribution in part (a) with $\beta=1$.
\end{Lemma}

We note that the \emph{independence} claimed in part (b) is non-trivial,
in particular, since the random variables $(\eta_m, \kappa_m, \rho_m)$ are
all defined in terms of the \emph{same} exponential random variables.
We shall later see a more general version of this result.

\proof Part (a) is quite trivial and we shall leave the proof to
the reader and focus on part (b). Note that for any fixed $x_0$, $y_0$ and $z_0$
all positive and for $X, Y, Z$ all independent exponential random variables,
we have
    \eqn{
    \lbeq{X+Y-est}
    \prob(X+Y\leq x_0/\sqrt{m})=\frac{x_0^2}{2m}+O(m^{-3/2}),
    }
and similar estimates hold for $\prob(X+ Z \leq y_0/\sqrt{m})$ and
$\prob(Y+ Z \leq z_0/\sqrt{m})$. Further, we make use of the fact that, for $m\rightarrow \infty$,
    \eqn{
    \prob\Big(X+Y\leq x_0/\sqrt{m} , X+ Z \leq y_0/\sqrt{m}\Big)=\Theta(m^{-3/2}),
    \label{eqn-three-exp}
    }
since $X+Y\leq x_0/\sqrt{m} , X+ Z \leq y_0/\sqrt{m}$ implies that $X,Y,Z$ are all of order $1/\sqrt{m}$.
Then, we rewrite
    \eqn{
    \prob\Big(\eta_m > x_0, \kappa_m > y_0, \rho_m> z_0\Big) = \prob\left(\sum_{i=1}^m I_i = 0, \sum_{i=1}^m J_i = 0, \sum_{i=1}^m L_i = 0\right),
    }
where $I_i = \ind_{\{X_i+ Y_i < x_0/\sqrt{m}\}}$ , $J_i =\ind_{\{X_i+ Z_i < y_0/\sqrt{m}\}}$
and $L_i = \ind_{\{Y_i+ Z_i < z_0/\sqrt{m}\}}$, where we write $\ind_{A}$ for the indicator of
the event $A$. This implies, in particular, that
    \eqan{
    \lbeq{union-rare-events}
    \prob(\eta_m > x_0, \kappa_m > y_0, \rho_m> z_0) &= \left(\prob(I_1 = 0 , J_1 = 0, L_1=0)\right)^m\\
    &=\left(1-\prob\big(\{I_1=1\}\cup \{J_1 =1\}\cup \{L_1=1\}\big)\right)^m\nn\\
    &=\left[1 - \left(\frac{x_0^2}{2{m}}+ \frac{y_0^2}{2{m}} +\frac{z_0^2}{2{m}}- \Theta(m^{-3/2})\right) \right]^m\nn\\
    &= {\mathrm e}^{-(x_0^2/2+y_0^2/2 + z_0^2/2)}(1+o(1)),\nn
    }
as $m\to \infty,$ where we use that
    \eqan{
    &\Big|\prob\big(\{I_1=1\}\cup \{J_1 =1\}\cup \{L_1=1\}\big)-\prob(I_1=1)-\prob(J_1=1)-\prob(L_1=1)\Big|\\
    &\qquad\leq \prob(I_1=J_1=1)+\prob(I_1=L_1=1)+\prob(J_1=L_1=1)
    =\Theta(m^{-3/2}).\nn
    }
This proves the result. \qed
\vskip0.5cm

\noindent
The next lemma generalizes the statement of Lemma \ref{lem-min-Gamma} in a substantial way:

\begin{Lemma}[Minima of Gamma random variables on the complete graph]
\label{lem-min-Gamma-gen}
Fix $k\geq 1$ and $n\geq k$. Let $\{E_{s,t}\}_{1\leq s<t\leq n}$ be an i.i.d.\ sequence of
exponential random variables with mean 1. For each $i\in [k]$, let $\NN_{i}\subseteq [n]\setminus [k]$
denote deterministic sets of indices.
Let $\NN_{ij}=\NN_i\cap \NN_j$, and assume that, for each $i,j\in [k]$,
    \eqn{
    \lbeq{conv-inters}
    |\NN_{ij}|/n\rightarrow \beta_{ij}>0.
    }
Let
    \eqn{
    \eta^{\sss(n)}_{ij}=\sqrt{n}\min_{s\in \NN_{ij}} (E_{i,s}+E_{s,j}).
    }
Then, for each $k$,
    \eqn{
    \{\eta_{ij}^{\sss (n)}\}_{1\leq i<j\leq k}
    \convd \{\eta_{ij}\}_{1\leq i<  j\leq k},
    }
where the random variables
$\{\eta_{ij}\}_{1\leq i<j\leq k}$ are independent
random variables with distribution
    \eqn{
    \prob(\eta_{ij} > x )
    \to \exp\left(-\beta_{ij} x^2/2\right).
    \lbeq{limit-weight}
    }
When $\NN_{i}$ denote \emph{random} sets of indices which are independent
of the exponential random variables, then the same result holds
when the convergence in \refeq{conv-inters} is replaced with
convergence in distribution where the limits $\beta_{ij}$ satisfy that
$\beta_{ij}>0$ holds a.s., and the limits $\{\eta_{ij}\}_{1\leq i<  j\leq k}$
are conditionally independent given $\{\beta_{ij}\}_{1\leq i<  j\leq k}$.
\end{Lemma}

\proof We follow the proof of Lemma \ref{lem-min-Gamma} as closely as possible.
For $i\in [k]$ and $s\in [n]\setminus [k]$, we define $X_{i,s}=E_{i,s},$ when $s\in \NN_i$,
and $X_{i,s}=+\infty,$ when $s\not\in \NN_i$.
Since the sets of indices $\{\NN_i\}_{i\in [k]}$ are independent
from the exponential random variables, the variables $\{X_{i,s}\}_{i\in [k], s\in [n]\setminus [k]}$
are, conditionally on $\{\NN_i\}_{i\in [k]}$,
\emph{independent} random variables. Then, since $\NN_{ij}=\NN_i\cap \NN_j$,
    \eqn{
    \eta^{\sss(n)}_{ij}=\sqrt{n}\min_{s\in \NN_{ij}} (E_{i,s}+E_{j,s})=
    \sqrt{n}\min_{s\in [n]\setminus [k]} (X_{i,s}+X_{j,s}).
    }
Let $\{x_{ij}\}_{1\leq i<j\leq k}$ be a vector with positive coordinates.
We note that
    \eqn{
    \prob(\eta^{\sss(n)}_{ij}>x_{ij},\forall i,j\in [k])
    =\prob\Big(\sum_{s\in [n]\setminus [k]} J_{ij,s}=0,\forall i,j\in [k]\Big),
    }
where $J_{ij,s}=\ind_{\{X_{i,s}+ X_{j,s} < x_{ij}/\sqrt{n}\}}$. We note that the random vectors
$\{J_{ij,s}\}_{s\in [n]\setminus [k]}$ are conditionally independent given
$\{\NN_i\}_{i\in [k]}$, so that
    \eqn{
    \prob(\eta^{\sss(n)}_{ij}>x_{ij},\forall i,j\in [k])
    =\prod_{s\in [n]\setminus [k]}\prob(J_{ij,s}=0,\forall i,j\in [k]).
    }
Now, note that $J_{ij,s}=0$ a.s.\ when $s\not\in \NN_{ij}$, while,
for $s\in \NN_{ij}$, we have, similarly to \refeq{X+Y-est},
    \eqn{
    \prob(J_{ij,s}=1)=\frac{x_{ij}^2}{2n}+O(n^{-3/2}).
    }
Therefore, we can summarize these two claims by
    \eqn{
    \prob(J_{ij,s}=1)= \ind_{\{s\in \NN_{ij}\}}\Big(\frac{x_{ij}^2}{2n}+\Theta(n^{-3/2})\Big).
    }
Similarly to the argument in \refeq{union-rare-events}, we have that
    \eqan{
    \prob(J_{ij,s}=0,\,\forall i,j\in [k])
    &=1-\sum_{1\leq i<j\leq k} \prob(J_{ij,s}=1)+\Theta(n^{-3/2})\nn\\
    &=\exp\Big\{-\sum_{1\leq i<j\leq k} \ind_{\{s\in \NN_{ij}\}}\Big(\frac{x_{ij}^2}{2n}+\Theta(n^{-3/2})\Big)\Big\}.
    }
We conclude that
    \eqan{
    \prob(\eta^{\sss(n)}_{ij}>x_{ij},\,\forall i,j\in [k])
    &=\prod_{s\in [n]\setminus [k]}\prob(J_{ij,s}=0\forall i,j\in [k])\\
    &=\exp\Big\{-\sum_{s\in [n]\setminus [k]}\sum_{1\leq i<j\leq k} \ind_{\{s\in \NN_{ij}\}}\Big(\frac{x_{ij}^2}{2n}+\Theta(n^{-3/2})\Big)\Big\}
    \nn\\
    &=\exp\{-\sum_{1\leq i<j\leq k} x_{ij}^2 \beta_{ij}/2\}(1+o(1)),\nn
    }
as required.
\qed
\vskip0.5cm

\noindent
We shall apply Lemma \ref{lem-min-Gamma-gen} to $\NN_i$ being the direct neighbors in $[n]\setminus [k]$
of vertex $i\in [k]$. Thus, by Lemma \ref{lem-min-Gamma-gen}, in order to prove the convergence of
the weights, it suffices to prove the convergence of the number of \emph{joint} neighbors of
the super vertices $i$ and $j$, simultaneously, for all $i,j\in [k]$. That is the content
of the following lemma:

\begin{Lemma}[Weak convergence of $\Ner_{ij}/n$]
\label{lem-fpipj-conv}
The random vector $\{\Ner_{ij}/n\}_{1\leq i<j\leq n},$
converges in distribution in the product topology to $\{f(P_i,P_j)\}_{1\leq i<j<\infty}$,
where $f(P_i,P_j)$ is defined in \refeq{fpipj-def}, and $\{P_i\}_{i\geq 1}$
has the Poisson-Dirichlet distribution.
\end{Lemma}

\proof We shall first prove that the random vector
$\{\Ner_{ij}/n-f\big(P_i^{\sss(n)},P_j^{\sss(n)}\big)\}_{1\leq i<j\leq n},$
converges in probability in the product topology to zero, where
$P_i^{\sss(n)}=D_{\sss (n+1-i:n)}/L_n$ is the normalized
$i^{\rm th}$ largest degree. For this, we note that
    \eqn{
    \Ner_{ij}=\sum_{s=1}^n I_{s}(i,j),
    }
where $I_{s}(i,j)$ is the indicator that $s\in [n]$ is a neighbor of both $i$ and $j$.
Now, weak convergence in the product topology is equivalent to the weak convergence of
$\{\Ner_{ij}/n\}_{1\leq i<j<K}$ for any $K\in \Zbold^+$ (see \cite[Theorem 4.29]{Kall02}).
For this, we shall use a second moment method.
We first note that $|\Ner_{ij}/n-\Ner_{\sss \leq b_n}(i,j)/n|
\leq \frac1n \sum_{s=1}^n \indic{D_s\geq b_n} \convp 0$,
where $b_n\rightarrow \infty$ and
    \eqn{
    \Ner_{\sss \leq b_n}(i,j)=\sum_{s=1}^n I_{s}(i,j)\indic{D_s\leq b_n}.
    }
Take $b_n=n$ and note that when $i,j\leq K$, the vertices $i$ and $j$ both have degree of order $n^{1/(\tau-1)}$ which is at least $n$ {\bf whp}.
Thus, the sum over $s$ in $N_{\sss \leq n}(i,j)$ involves different vertices than $i$ and $j$. Next,
we note that
    \eqan{
    \expec_n[\Ner_{\sss \leq n}(i,j)/n]&=\frac 1n \sum_{s=1}^n \indic{D_s\leq n} \prob_n(I_{s}(i,j)=1)\nn\\
    &=\frac 1n \sum_{s=1}^n \indic{D_s\leq n}[1-(1-P_i^{\sss(n)})^{D_s}-(1-P_j^{\sss(n)})^{D_s}
    +(1-P_i^{\sss(n)}-P_j^{\sss(n)})^{D_s}]+\op(1),
    }
in a similar way as in \refeq{fpipj-def}. By dominated convergence, we have that,
for every $s\in [0,1]$,
    \eqn{
    \frac 1n \sum_{s=1}^n \indic{D_s\leq n}(1-s)^{D_s}\convas \expec[(1-s)^D],
    }
which implies that
    \eqn{
    \expec_n[\Ner_{\sss \leq n}(i,j)/n]-f\big(P_i^{\sss(n)},P_j^{\sss(n)}\big)\convp 0.
    }
Further, the indicators $\{I_{s}(i,j)\}_{s=1}^n$ are close to independent, so that
${\rm Var}_n\big(\Ner_{\sss \leq n}(i,j)/n\big)=\op(1)$, where ${\rm Var}_n$
denotes the variance w.r.t.\ $\prob_n$.
The weak convergence claimed in Lemma \ref{lem-fpipj-conv} follows
directly from the above results, as well as the weak convergence of the
order statistics in \refeq{EVTexpsrep} and the continuity of $(s,t)\mapsto f(s,t)$.
\qed
\medskip

\noindent
The following corollary completes the proof of the convergence of the rescaled minimal
weight two-edge paths in $\GGer_n$:

\begin{Corollary}[Conditional independence of weights]
\label{cor-cond-ind-weights}
Let $l_{ij}^{\sss (n)} = \sqrt{n} w_{ij}^{\sss (2)}$, where $w_{ij}^{\sss (2)}$ is the minimal
weight of all two-edge paths between the vertices $i$ and $j$ (with $w_{ij}^{\sss (2)} = \infty$
if they are not connected by a two-edge path). Fix $k\geq 1$. Then,
    \eqn{\left(\{l_{ij}^{\sss (n)}\}_{1\leq i<  j\leq k},
    \{{D_i}/{L_n}\}_{1\leq i\leq n} \right) \convd \left(\{l_{ij}\}_{1\leq i< j\leq k} ,
    \{P_i\}_{i\geq 1}\right),}
where, given $\{P_i\}_{i\geq 1}$
the random variables  $\{l_{ij}\}_{1\leq i<j\leq k}$ are conditionally independent
with distribution
    \eqn{
    \prob(l_{ij}>x )
    \to \exp\left(-f(P_i,P_j) x^2/2\right).
    \lbeq{eq:lim-lij-KKer}
    }
\end{Corollary}

\proof The convergence of $ \{{D_m}/{L_n}\}_{1\leq m\leq n}$ follows from Section \ref{sec-PDdistr}.
Then we apply Lemma \ref{lem-min-Gamma-gen}. We let $\NN_i$ denote the set of neighbors in
$[n]\setminus [k]$ of the super vertex $i\in [k]$. Then, $|\NN_{ij}|
=|\NN_i\cap \NN_j|=Ner_{ij}$, so that \refeq{conv-inters} is equivalent to
the convergence in distribution of $\Ner_{ij}/n$. The latter is proved in Lemma
\ref{lem-fpipj-conv}, with $\beta_{ij}=f(P_i,P_j)$. Since $P_i>0$ a.s.\ for each $i\in [k]$,
we obtain that $\beta_{ij}>0$ a.s.\ for all $i,j\in [k]$. Therefore, Lemma \ref{lem-min-Gamma-gen}
applies, and completes the proof of the claim.
\qed
\medskip

\noindent
Now we are ready to prove Proposition \ref{prop-weak-con-FPP-erased}:\\
{\it Proof of Proposition \ref{prop-weak-con-FPP-erased}.}
By Corollary \ref{cor-cond-ind-weights}, we see that the weights in the
FPP problem $\KK_n^k$ converge in distribution to the weights in the FPP
on $(\KKer)^k$. Since the weights $W_{ij}^{\sss(n)}$ of the minimal two-edge
paths between $i,j\in [k]$
are continuous functions of the weights $\{l^{\sss(n)}_{ij}\}_{1\leq i<j\leq k}$,
it follows that  $\{W_{ij}^{\sss(n)}\}_{1\leq i<j\leq k}$ converges in
distribution to $\{W_{ij}\}_{1\leq i<j\leq k}$.
Since the weights are continuous random variables, this also implies that
the hopcounts $\{H^{\sss(n)}_{ij}\}_{1\leq i<j\leq k}$ in $\KK_n^k$
converge in distribution to the hopcounts $\{H_{ij}\}_{1\leq i<j\leq k}$
in $(\KKer)^k$. This proves Proposition \ref{prop-weak-con-FPP-erased}.
\qed

\subsection{Coupling of the minimal edges from uniform vertices: Proof of Proposition \ref{prop-coupling+tightness-weight}}
\label{sec-pf-prop6.2}
In this section, we prove Proposition \ref{prop-coupling+tightness-weight}.
We start by noticing that the vertices $\UnVer_i,\, i=1,2,$ are, {\bf whp}, only attached to super
vertices. Let $I^{\sss(n)}$ and $J^{\sss(n)}$ denote the vertices to which $\UnVer_i,\, i=1,2,$ are connected
and of which the edge weights are minimal. Then, by the discussion below \eqref{eq:lij},
$(I^{\sss(n)},J^{\sss(n)})$ converges in distribution to the random vector
$(I,J)$ having the distribution specified right before Theorem \ref{thm-erased}, and
where the two components are  conditionally independent, given $\{P_i\}_{i\geq 1}$.

Further, denote the weight of the edges attaching $(\UnVer_1,\UnVer_2)$ to $(I^{\sss(n)},J^{\sss(n)})$
by $(V^{\sss(n)}_1, V^{\sss(n)}_2)$. Then, $(V^{\sss(n)}_1, V^{\sss(n)}_2)\convd (\Ver_1,\Ver_2)$
defined in Theorem \ref{thm-erased}. This in particular proves \refeq{Wn-weak-conv-simple}
since the weight between any two super vertices is $\op(1)$.
Further, since $(I^{\sss(n)},J^{\sss(n)})$ are discrete random variables that weakly converge
to $(I,J)$, we can couple $(I^{\sss(n)},J^{\sss(n)})$ and $(I,J)$ in such a way that
\refeq{coupling-IJ} holds.

Let $(D^{{\rm er}\sss(n)}_{\UnVer_1}, D^{{\rm er}\sss(n)}_{\UnVer_2})$ denote the erased degrees
of the vertices $(\UnVer_1,\UnVer_2)$ in $\GGer$. The following lemma states
that these erased degrees converge in distribution:

\begin{Lemma}[Convergence in distribution of erased degrees]
\label{lem-conv-d-erased-degrees}
Under the conditions of Theorem \ref{thm-erased}, as $n\rightarrow \infty$,
    \eqn{
    \lbeq{conv-D-er-2}
    (D^{{\rm er}\sss(n)}_{\UnVer_1}, D^{{\rm er}\sss(n)}_{\UnVer_2})
    \convd (\Der_1,\Der_2),
    }
which are two copies of the random variable $\Der$
described right before Theorem \ref{thm-erased}, and
which are conditionally independent given $\{P_i\}_{i\geq 1}$.
\end{Lemma}

\proof We note that the degrees before erasure, i.e., $(D_{\UnVer_1}, D_{\UnVer_2})$,
are i.i.d.\ copies of the distribution $D$ with distribution function $F$, so that,
in particular, $(D_{\UnVer_1}, D_{\UnVer_2})$ are bounded by $K$ {\bf whp} for any $K$
sufficiently large. We next investigate the effect of erasure. We condition on
$\{P_i^{\sss(n)}\}_{i=1}^{m_n}$, the rescaled $m_n$ largest degrees, and note that, by
\refeq{EVTexpsrep}, $\{P_i^{\sss(n)}\}_{i=1}^{m_n}=\{D_i/L_n\}_{i=1}^{m_n}$
converges in distribution to $\{P_i\}_{i\geq 1}$. We let $m_n\rightarrow \infty$
arbitrarily slowly, and note that, {\bf whp}, the $(D_{\UnVer_1}, D_{\UnVer_2})$
half-edges incident to the vertices  $(\UnVer_1,\UnVer_2)$, are exclusively connected
to vertices in $[m_n]$. The convergence in \refeq{conv-D-er-2} follows when
    \eqan{
    \lbeq{aim-conv}
    &\prob\Big((D^{{\rm er}\sss(n)}_{\UnVer_1}, D^{{\rm er}\sss(n)}_{\UnVer_2})=(k_1,k_2)\mid
    \{P_i^{\sss(n)}\}_{i=1}^{m_n}, (D_{\UnVer_1}, D_{\UnVer_2})=(j_1,j_2)\Big)\\
    &\qquad\qquad =G_{k_1,j_1}(\{P_i^{\sss(n)}\}_{i=1}^{m_n})G_{k_2,j_2}(\{P_i^{\sss(n)}\}_{i=1}^{m_n})\nn
    +\op(1),
    }
for an appropriate function $G_{k,j}\colon {\mathbb R}_+^{\N} \to [0,1]$, which, for every $k,j$,
is \emph{continuous} in the product topology. (By convention, for a vector with finitely many coordinates $\{x_i\}_{i=1}^m$, we let $G_{k_1,j_1}(\{x_i\}_{i=1}^m)=G_{k_1,j_1}(\{x_i\}_{i=1}^{\infty})$, where
$x_i=0$ for $i>m$.)

Indeed, from \refeq{aim-conv}, it follows that, by dominated
convergence,
    \eqan{
    \prob\Big((D^{{\rm er}\sss(n)}_{\UnVer_1}, D^{{\rm er}\sss(n)}_{\UnVer_2})=(k_1,k_2)\Big)
    &=\expec\Big[\prob\Big((D^{{\rm er}\sss(n)}_{\UnVer_1}, D^{{\rm er}\sss(n)}_{\UnVer_2})=(k_1,k_2)\mid
    \{P_i^{\sss(n)}\}_{i=1}^{m_n}, (D_{\UnVer_1}, D_{\UnVer_2})\Big)\Big]\nn\\
    &=\expec\Big[G_{k_1,D_1}(\{P_i^{\sss(n)}\}_{i=1}^{m_n})G_{k_2,D_2}(\{P_i^{\sss(n)}\}_{i=1}^{m_n})\Big]
    +o(1)\nn\\
    &\rightarrow \expec\Big[G_{k_1,D_1}(\{P_i\}_{i\geq 1})G_{k_2,D_2}(\{P_i\}_{i\geq 1})\Big],
    }
where the last convergence follows from weak convergence of
$\{P_i^{\sss(n)}\}_{i=1}^{m_n}$ and the assumed continuity of $G$.
The above convergence, in turn, is equivalent to \refeq{conv-D-er-2},
when $G_{k,j}(\{P_i\}_{i\geq 1})$ denotes the probability that
$k$ distinct cells are chosen in a multinomial experiment with $j$
independent trials where, at each trial,  we choose cell $i$ with
probability $P_i$. It is not hard to see that, for each $k,j$,
$G_{k,j}$ is indeed a continuous function in the product topology.

To see \refeq{aim-conv}, we note that, conditionally on
$\{P_i^{\sss(n)}\}_{i=1}^{m_n}$, the vertices to which the $D_{\UnVer_i}=j_i$
stubs attach are close to independent, so that it suffices to prove that
    \eqn{
    \lbeq{aim-conv-rep}
    \prob\big(D^{{\rm er}\sss(n)}_{\UnVer_1}=k_1\mid
    \{P_i^{\sss(n)}\}_{i=1}^{m_n}, D_{\UnVer_1}=j_1\big)=G_{k_1,j_1}(\{P_i^{\sss(n)}\}_{i=1}^{m_n})
    +\op(1).
    }
The latter follows, since, again conditionally on
$\{P_i^{\sss(n)}\}_{i=1}^{m_n}$, each stub chooses to connect to
vertex $i$ with probability $D_i/L_n=P_i^{\sss(n)}$, and the different
stubs choose close to independently. This completes the proof of
Lemma \ref{lem-conv-d-erased-degrees}.
\qed

By Lemma \ref{lem-conv-d-erased-degrees}, we can also couple
$(D^{{\tt er}\sss(n)}_{\UnVer_1}, D^{{\tt er}\sss(n)}_{\UnVer_2})$
to $(\Der_1,\Der_2)$ in such a way that
    \eqn{
    \lbeq{coupling-Ders}
    \prob\big((D^{{\tt er}\sss(n)}_{\UnVer_1}, D^{{\tt er}\sss(n)}_{\UnVer_2})\neq (\Der_1,\Der_2)\big)=o(1).
    }
Now, $(V^{\sss(n)}_1, V^{\sss(n)}_2)$ is equal in distribution to
$(E_1/D^{{\tt er}\sss(n)}_{\UnVer_1}, E_2/D^{{\tt er}\sss(n)}_{\UnVer_2})$, where
$(E_1,E_2)$ are two independent exponential random variables with mean 1.
Let $V_i=\Ver_i=E_i/\Der_i$, where we use the \emph{same} exponential
random variables. Then $(V_1,V_2)$ has the right distribution, and
the above coupling also provides a coupling of
$(V^{\sss(n)}_1, V^{\sss(n)}_2)$ to $(V_1,V_2)$ such that
\refeq{coupling-Vs} holds.

By the above couplings, we have that
$\sqrt{n}\big(W_n-(V^{\sss(n)}_1+V^{\sss(n)}_2)\big)=
\sqrt{n}\big(W_n-(V_1+V_2)\big)$ {\bf whp}.
By construction, $\sqrt{n}\big(W_n-(V^{\sss(n)}_1+V^{\sss(n)}_2)\big)\geq 0$ a.s.,
so that also, {\bf whp}, $\sqrt{n}\big(W_n-(V_1+V_2)\big)\geq 0$. Further,
$\sqrt{n}\big(W_n-(V^{\sss(n)}_1+V^{\sss(n)}_2)\big)\leq l^{\sss(n)}_{I^{\sss(n)},J^{\sss(n)}},$
which is the weight of the minimal two-edge path between the super vertices
$I^{\sss(n)}$ and $J^{\sss(n)}$. Now, by \refeq{coupling-IJ},
$(I^{\sss(n)},J^{\sss(n)})=(I,J)$ {\bf whp}. Thus, {\bf whp},
$l^{\sss(n)}_{I^{\sss(n)},J^{\sss(n)}}=l^{\sss(n)}_{I,J}$, which, by Proposition
\ref{prop-weak-con-FPP-erased}, converges in distribution to $l_{IJ}$, which is
a finite random variable. As a result, $l^{\sss(n)}_{I^{\sss(n)},J^{\sss(n)}}$ is a
tight sequence of random variables, and, therefore, also
$\sqrt{n}\big(W_n-(V^{\sss(n)}_1+V^{\sss(n)}_2)\big)$ is.
This completes the proof of Proposition \ref{prop-coupling+tightness-weight}.
\qed

\subsection{Tightness of FPP problem and evenness of hopcount: Proof of Proposition \ref{prop-hopcount-conv}}
\label{sec-pf-prop6.3}
In this section, we prove that the only possible minimal weight paths between the
super vertices are two-edge paths. All other paths are much too costly to be used.
We start by stating and proving a technical lemma
about expectations of degrees conditioned to be at most $x$.
It is here that we make use of the condition in \refeq{assumption-F-erased}:

\begin{Lemma}[Bounds on restricted moments of $D$]
\label{lem-rest-mom-D}
Let $D$ be a random variable with distribution function $F$ satisfying
\refeq{assumption-F-erased} for some $\tau\in (1,2).$ Then,
there exists a constant $C$ such that, for every $x\geq 1$,
    \eqn{
    \lbeq{rest-mom-D}
    \expec[D\indic{D\leq x}]\leq C x^{2-\tau},
    \quad
    \expec[D^{\tau-1}\indic{D\leq x}]\leq C\log{x},
    \quad
    \expec[D^{\tau}\indic{D\leq x}]\leq Cx,
    \quad
    \expec[D^{2(\tau-1)}\indic{D\leq x}]\leq Cx^{\tau-1}.
    }
\end{Lemma}

\proof We note that,
for every $a>0$, \rr{using partial integration},
    \eqn{
    \lbeq{pi}
    \expec[D^a\indic{D\leq x}]=
    -\int_{(0,x]}y^a\,d(1-F(y))
    \leq a \int_0^x y^{a-1}[1-F(y)]dy\leq c_2 a \int_0^x y^{a-\tau}dy .
    }
The proof is completed by considering the four cases separately and computing in each case the integral on the right-hand side of \refeq{pi}.
\qed
\medskip

\noindent
The following lemma shows that paths of an {\it odd} length are unlikely:

\begin{Lemma}[Shortest-weight paths on super vertices are of even length]
\label{lem-two-paths-odd}
Let the distribution function $F$ of the degrees of the
CM satisfy \refeq{assumption-F-erased}. Let ${\cal B}^{\sss(n)}$ be the event
that there exists a path between two super vertices consisting of all normal
vertices and having an odd number of edges and of total
weight $w_n/\sqrt{n}$. Then, for some constant $C$,
    \eqn{
    \lbeq{two-path-prob}
    \prob({\cal B}^{\sss(n)})\leq \frac{\eps_n^{-2(\tau-1)}}{\sqrt{n\log{n}}}{\mathrm e}^{C w_n \sqrt{\log{n}}}.
    }
\end{Lemma}

\proof We will show that the probability that there exists
a path between two super vertices consisting of all normal
vertices and having an odd number of edges and of total
weight $w_n/\sqrt{n}$ is small. For this, we shall use the
first moment method and show that the expected number of
such paths goes to $0$ as $n\to\infty$. Fix two super vertices which will be the end
points of the path and an {\it even} number $m\geq 0$ of normal vertices
with indices $i_1, i_2, \ldots i_m$. Note that when a path \rr{between} two super vertices
consists of an even number of vertices, then the path has an odd number of
edges.

Let $\CBnm$ be the event that there
exists a path between two super vertices consisting of exactly $m$
intermediate normal vertices with total weight $w_n/\sqrt{n}$.
We start by investigating the case $m=0$, so that the super vertices
are directly connected.
Note that $|\su_n|=\Op(\expec[|\su_n|])$, by concentration, and that
    $$
    \expec[|\su_n|]=n\prob(D_1>\vep_n n^{1/(\tau-1)})=O(\eps_n^{-(\tau-1)}),
    $$
Hence, there are $\Op(\eps_n^{-(\tau-1)})$ super vertices
and thus $\Op(\eps_n^{-2(\tau-1)})$ edges between them. The probability that any one of
them is smaller than $w_n/\sqrt{n}$ is of order $\vep_n^{-2(\tau-1)}w_n/\sqrt{n}$, and
it follows that $\prob({\cal B}_0^{\sss(n)})\leq \vep_n^{-2(\tau-1)}w_n/\sqrt{n}$.

Let $\Mnm$ be the total number of paths connecting two specific super vertices
and which are such that the total weight on the paths is at most $w_n/\sqrt{n}$,
so that
    \eqn{
    \prob(\CBnm)=\prob(\Mnm \geq 1)\leq \expec[\Mnm].
    }
In the following argument, for convenience, we let $\{D_i\}_{i=1}^n$ denote the i.i.d.\ vector of
degrees (i.e., below $D_i$ is not the $i^{\rm th}$ largest degree, but rather
a copy of the random variable $D\sim F$ independently of the other degrees.)

Let $\veciota=(i_1, i_2,\ldots, i_m)$, and denote
by $p_{m,n}(\veciota)$ the probability that the $m$ vertices $i_1, i_2, \ldots, i_m$
are normal and are such that there is an edge between $i_{s}$ and $i_{s+1},$ for $s=1, \ldots, m-1$.
Further, note that with $S_{m+1} = \sum_{i=1}^{m+1} E_i$, where $E_i$ are
independent exponential random variables \rr{with mean 1}, we have, for any $u\in [0,1]$,
    \eqn{
    \lbeq{vijfachttien}
    \prob(S_{m+1}\leq u) =\int_0^u \frac{x^m {\mathrm e}^{-x}}{m!} \leq  \frac{u^{m+1}}{(m+1)!}.
    }
Together with the fact that there are \rr{$\Op(\eps_n^{-(\tau-1)})$} super vertices, this implies that
    \eqn{
    \label{eqn:bmn}
    \prob(\CBnm)\leq
    \expec[\Mnm]\leq \frac{C\eps_n^{-2(\tau-1)} w_n^{m+1}}{(m+1)! n^{(m+1)/2}}
    \sum_{\veciota} p_{m,n}(\veciota),
    }
since \refeq{vijfachttien} implies that the probability that the sum
of $ m+1$  exponentially distributed  r.v.'s is smaller than
$u_n=w_n/\sqrt{n}$ is at most $ u_n^{m+1}/(m+1)!$.

By the construction of the CM, we have
    \eqn{
    \lbeq{pmn-bound}
    p_{m,n}(\veciota)\leq \expec\left[\prod_{j=1}^{m-1}\left(\frac{D_{i_j} D_{i_{j+1}}}{L_n-2j+1} \wedge 1\right)\indicwo{\FF_m}\right]
    \leq \expec\left[\prod_{j=1}^{m-1}\left(\frac{D_{i_j} D_{i_{j+1}}}{L_n} \wedge 1\right)\indicwo{\FF_m}\right](1+o(1)),
    }
where $\FF_m$ is the event that $D_{i_j} < \eps_n n^{1/(\tau-1)}$ for all $1\leq j\leq m$.
We shall prove by induction that, for every $\veciota$, and for $m$ even,
    \eqn{
    \lbeq{IH-pmn}
    p_{m,n}(\veciota)   \leq
    \frac{(C \log{n})^{m/2}}{n^{m/2}}.
    }
We shall initiate \refeq{IH-pmn} by verifying it for $m=2$ directly, and
then advance the induction by relating $p_{m,n}$ to $p_{m-2,n}$.

We start by investigating expectations as in \refeq{pmn-bound} iteratively.
First, conditionally on $D_{i_{m-1}}$,
note that
\rr{
    \eqan{
    \lbeq{expectation-one}
    \expec\Big[\frac{D_{i_{m-1}} D_{i_m}}{L_n}\wedge 1\Big|D_{i_{m-1}}\Big]
    &= \prob\Big(D_{i_m}>\frac{L_n}{D_{i_{m-1}}}\big|D_{i_{m-1}}\Big)
    +D_{i_{m-1}} \expec\big[\frac{D_{i_m}}{L_n} \indic{D_{i_m}\leq L_n/D_{i_{m-1}}}\big|D_{i_{m-1}}\big]\nn\\
    }
    }
Furthermore,
    \eqn{
    \prob\Big(D_{i_m}>\frac{L_n}{D_{i_{m-1}}}\big|D_{i_{m-1}}\Big)
    \le
     c_2(D_{i_{m-1}})^{\tau-1}\expec\Big[(L_n)^{1-\tau}\big|D_{i_{m-1}}\Big].
    }
In a similar way, we obtain using the first bound in Lemma \ref{lem-rest-mom-D}
together with the fact that $\{D_j\}_{j=1}^n$ is an i.i.d.\ sequence,
that
    \eqan{
    D_{i_{m-1}} \expec\big[\frac{D_{i_m}}{L_n} \indic{D_{i_m}\leq L_n/D_{i_{m-1}}}\big|D_{i_{m-1}}\big]
    &\leq C D_{i_{m-1}}\expec\Big[L_n^{-1}(L_n/D_{i_{m-1}})^{2-\tau}\big|D_{i_{m-1}}\Big]\nn\\
    &=C (D_{i_{m-1}})^{\tau-1}\expec\Big[(L_n)^{1-\tau}\big|D_{i_{m-1}}\Big],
    }
where we reach an equal upper bound as above. Thus,
    \eqan{
    \lbeq{expectation-one-rep}
    \expec\Big[\frac{D_{i_{m-1}} D_{i_m}}{L_n}\wedge 1\Big|D_{i_{m-1}}\Big]
    &\leq C (D_{i_{m-1}})^{\tau-1}\expec\Big[(L_n)^{1-\tau}\big|D_{i_{m-1}}\Big].
    }
Now, \cite[Lemma 4.1(b)]{DeiEskHofHoo09} implies that
$\expec[(L_n)^{1-\tau}|D_{i_{m-1}}]\le \expec[(L_n-D_{i_{m-1}})^{-(\tau-1)}]\leq c/n$, a.s.
so that
    \eqn{
    \expec\Big[\prob\Big(D_{i_m}>\frac{L_n}{D_{i_{m-1}}}\indic{D_{i_{m-1}}\leq \vep_n n^{1/(\tau-1)}}\big|D_{i_{m-1}}\Big)\Big]
    \le C \log{n}/n ,
    }
where, in the inequality, we have used the second
inequality in Lemma \ref{lem-rest-mom-D} together with the fact that $\{D_j\}_{j=1}^n$ is an i.i.d.\ sequence.
The second term on the right-hand side of \refeq{expectation-one} can be treated similarly, and yields the same upper
bound. Putting the two bounds together we arrive at
    \eqan{
    \lbeq{IH-m2}
    p_{2,n}(i_1,i_2)&=\expec\left[\left(\frac{D_{i_1} D_{i_2}}{L_n} \wedge 1\right)\indicwo{\FF_2}\right]
    \leq C\log{n}/n.
    }
which is \refeq{IH-pmn} for $m=2$.

To advance the induction, we need to extend \refeq{expectation-one}.
Indeed, we use \refeq{expectation-one-rep} to compute that
    \eqan{
    \lbeq{expectation-two-rep}
    &\expec\left[\left(\frac{D_{i_{m-2}} D_{i_{m-1}}}{L_n} \wedge 1\right)\cdot\left(\frac{D_{i_{m-1}} D_{i_m}}{L_n}\wedge 1\right)\indicwo{\FF_m}
    \big|D_{i_{m-2}}\right]\\
    &\qquad \leq C^\prime \expec\left[\left(\frac{D_{i_{m-2}} D_{i_{m-1}}}{L_n} \wedge 1\right) \left(\frac{D_{i_{m-1}}}{L_n}\right)^{\tau-1}\indicwo{\FF_m}
    \big|D_{i_{m-2}}\right]\nn\\
    &\qquad =C^\prime \expec\left[\left(\frac{D_{i_{m-1}}}{L_n}\right)^{\tau-1} \indic{D_{i_{m-1}} > \frac{L_n}{D_{i_{m-2}}}}
    \indicwo{\FF_{m-1}}\big|D_{i_{m-2}}\right]\nn\\
    &\qquad\quad+ C^\prime \expec\left[\left(\frac{D_{i_{m-1}}}{L_n}\right)^{\tau} D_{i_{m-2}}\indic{D_{i_{m-1}} < \frac{L_n}{D_{i_{m-2}}}}
    \indicwo{\FF_{m-1}}\big|D_{i_{m-2}}\right]\nn,
    }
Now using Lemma \ref{lem-rest-mom-D} together with the fact that $\{D_j\}_{j=1}^n$ is an i.i.d.\ sequence, and simplifying, we obtain
the following to hold {\it almost surely},
    \eqn{
    \lbeq{expec-two-two-ways}
    \expec\left[\left(\frac{D_{i_{m-2}} D_{i_m-1}}{L_n} \wedge 1\right)
    \cdot\left(\frac{D_{i_{m-1}} D_{i_m}}{L_n}\wedge 1\right)\indicwo{\FF_m}\big|
    D_{i_{m-2}},D_{i_m}\right]
    \leq C\log{n}/n.
    }
This shows that $p_{m,n}\le (C\log{n}/n) p_{m-2,n}$, and hence proves \refeq{IH-pmn}.

Using this estimate in \eqref{eqn:bmn}, and summing over all even $m$, using the notation that $m=2\Zbold^+$, shows that
    \begin{eqnarray}
    \prob({\cal B}^{\sss(n)})&\leq& \sum_{m=2\Zbold^+}
    \prob(\CBnm)
    \leq \sum_{m=2\Zbold^+}
    \frac{C\eps_n^{-2(\tau-1)} w_n^{m+1}}{(m+1)! n^{(m+1)/2}}
    \sum_{\veciota}   p_{m,n}(\veciota)\nn\\
    & \leq&\sum_{m=2\Zbold^+}
    \frac{C\eps_n^{-2(\tau-1)} w_n^{m+1}}{(m+1)! n^{(m+1)/2}}
    (C n\log{n})^{m/2}=    C \eps_n^{-2(\tau-1)} n^{-1/2}
    \sum_{k=1}^{\infty} \frac{w_n^{2k+1}(C \log{n})^{k}}{(2k+1)!}\nn\\
    &\leq& C\frac{\eps_n^{-2(\tau-1)}}{n^{1/2}\sqrt{\log{n}}}{\mathrm e}^{C w_n \sqrt{\log{n}}}.
    \end{eqnarray}
\qed
\medskip

\noindent Lemma \ref{lem-two-paths-odd} shows that \rr{with the
correct choice of $\vep_n$, we find that} $\prob(H_n\not\in
2\Zbold^+)=o(1)$, and to prove Theorem \ref{thm-erased}, we shall
show that the shortest-weight paths between any two specific super
vertices alternate between super vertices and normal vertices. We will prove this statement,
in Lemma \ref{lem-two-paths-tightness} below, by showing that the probability that a  vertex with index at least
$K$ is used at an even place, is for $K$ large, quite small. This shows in particular
that, {\bf whp}, at all even places we have super vertices.  In the following
lemma, we collect the properties of the degrees and erased degrees
that we shall make use of in the sequel. In its statement, we
define
    \eqn{
    {\cal G}^{\sss(n)}={\cal G}^{\sss(n)}_1\cap {\cal G}^{\sss(n)}_2\cap {\cal G}^{\sss(n)}_3,
    }
where, for $a\in (0,1)$ and $C, \Cer>0$, we let
    \eqan{
    {\cal G}^{\sss(n)}_1&=\big\{L_n n^{-1/(\tau-1)}\in [a, a^{-1}]\big\},\\
    {\cal G}^{\sss(n)}_2&=\Big\{C^{-1}(n/i)^{1/(\tau-1)}\leq D_{\sss(n+1-i:n)} \leq C(n/i)^{1/(\tau-1)}, \forall i\in [n]\Big\},\\
    {\cal G}^{\sss(n)}_3&=\Big\{\Der_i\leq \Cer (n/i), \forall i\in [n]\Big\}.
    }
The event ${\cal G}^{\sss(n)}$ is the \emph{good event} that we shall work with.
We shall first show that, if we take $a>0$ sufficiently small and $C,\Cer$ sufficiently
large, then $\prob({\cal G}^{\sss(n)})$ is close to 1:

\begin{Lemma}[The good event has high probability]
\label{lem-degrees-in-erased}
For every $\vep>0$, there exist $a>0$ sufficiently small and
$C,\Cer$ sufficiently large such that
    \eqn{
    \prob({\cal G}^{\sss(n)})\geq 1-\vep.
    }
\end{Lemma}

\proof We split
    \eqn{
    \prob\big(({\cal G}^{\sss(n)})^c\big)
    =\prob\big(({\cal G}^{\sss(n)}_1)^c\big)
    +\prob\big( ({\cal G}^{\sss(n)}_2)^c\big)
    +\prob\big({\cal G}^{\sss(n)}_1\cap {\cal G}^{\sss(n)}_2\cap ({\cal G}^{\sss(n)}_3)^c\big),
    }
and bound each term separately. We can make $\prob\big(({\cal G}^{\sss(n)}_1)^c\big)\leq \vep/3$
by choosing $a>0$ sufficiently small by the weak convergence in \refeq{mean-PDD}.

To bound $\prob\big({\cal G}^{\sss(n)}_1\cap ({\cal G}^{\sss(n)}_2)^c\big)$,
we note that $D_{\sss(n+1-l:n)}>C(n/l)^{1/(\tau-1)}$ is equivalent to the statement that
the number of values $i$ such that $D_i > C (n/l)^{1/(\tau-1)}$ is at least
$l$. Since $\{D_i\}_{i=1}^n$ is an i.i.d.\
sequence, this number has a Binomial distribution with parameters $n$ and
success probability
    \eqn{
    q_{l,n}=[1-F(C(n/l)^{1/(\tau-1)})]\leq c_2 C^{-(\tau-1)} l/n,
    }
by \refeq{assumption-F-erased}. When the mean of this binomial,
which is $c_2 C^{-(\tau-1)} l$ is much smaller than $l$, which is equivalent to
$C>0$ being large, the probability
that this binomial exceeds $l$ is exponentially small in $l$:
    \eqn{
    \prob(D_{\sss(n+1-l:n)}>C(n/l)^{1/(\tau-1)})\leq {\mathrm e}^{-I(C) l},
    }
where $I(C)\rightarrow \infty$ when $C\rightarrow \infty$.
Thus, by taking $C$ sufficiently large, we can make the probability that
there exists an $l$ for which $D_{\sss(n+1-l:n)}>C(n/l)^{1/(\tau-1)}$ small.
In more detail,
    \eqn{
    \prob\Big(\exists l: D_{\sss(n+1-l:n)}>C(n/l)^{1/(\tau-1)}\Big)
    \leq \sum_{l\in [n]} \prob(D_{\sss(n+1-l:n)}>C(n/l)^{1/(\tau-1)})\leq
    \sum_{l\in [n]}{\mathrm e}^{-I(C) l} \leq \vep/3,
    }
when we make $C>0$ sufficiently large. In a similar way, we can show that
the probability that there exists $l$ such that $D_{\sss(n+1-l:n)}\leq C^{-1}(n/l)^{1/(\tau-1)}$
is small when $C>0$ is large.

In order to bound $\prob\big({\cal G}^{\sss(n)}_1\cap {\cal G}^{\sss(n)}_2\cap ({\cal G}^{\sss(n)}_3)^c\big)$
we need to investigate the random variable $\Der_i$. We claim that
there exists a $R=R(a,C,\Cer)$ with $R(a,C,\Cer)\rightarrow \infty$ as
$\Cer\rightarrow \infty$ for each fixed
$a,C>0$, such that
    \eqn{
    \lbeq{bd-Der-cond}
    \prob\big(\Der_i\geq \Cer j^{\tau-1}|D_i=j, {\cal G}^{\sss(n)}_1\cap {\cal G}^{\sss(n)}_2)\leq {\mathrm e}^{-R j^{\tau-1}}.
    }

Fix $\Cer>0$. In order for $\Der_i\geq \Cer j^{\tau-1}$ to occur,
we must have that at least $\Cer j^{\tau-1}/2$ of the neighbors of vertex $i$
have index at least $\Cer j^{\tau-1}/2$, where we recall that vertex $i$ is
such that $D_i=D_{\sss(n+1-i:n)}$ is the $i^{\rm th}$ largest degree.
The $j$ neighbors of vertex $i$ are close to being independent, and
the probability that any of them connects to a vertex with index at
least $k$ is, conditionally on the degrees $\{D_i\}_{i=1}^n$, equal to
    \eqn{
    \sum_{l\geq k} D_{\sss(n+1-l:n)}/L_n.
    }
When ${\cal G}^{\sss(n)}_1\cap {\cal G}^{\sss(n)}_2$ holds, then
$D_{\sss(n+1-l:n)}/L_n\leq (C/a) l^{-1/(\tau-1)}$, so
that
    \eqn{
    \sum_{l\geq k} D_{\sss(n+1-l:n)}/L_n\leq c(C/a)k^{-(2-\tau)/(\tau-1)}.
    }
As a result, we can bound the number of neighbors of vertex $i$ by a binomial
random variable with $p=c'k^{-(2-\tau)/(\tau-1)}$, where $k=\Cer j^{\tau-1}/2$, i.e.,
    \eqn{
    \prob\big(\Der_i\geq \Cer j^{\tau-1}|D_i=j,{\cal G}^{\sss(n)}_1\cap {\cal G}^{\sss(n)}_2)
    \leq \prob\big({\rm Bin}(j, c(C/a) j^{-(2-\tau)})\geq \Cer j^{\tau-1}/2\big).
    }
Next, we note that the mean of the above binomial random variable is
given by $c(C/a)j^{1-(2-\tau)}=c(C/a) j^{\tau-1}$. A concentration result for binomial random variables \cite{Jans02},
yields that for $C>0$ sufficiently large,
    \eqn{
    \prob\big(\Der_i\geq \Cer j^{\tau-1}|D_i=j, {\cal G}^{\sss(n)}_1\cap {\cal G}^{\sss(n)}_2)
    \leq {\mathrm e}^{-R j^{\tau-1}}.
    }
This proves \refeq{bd-Der-cond}.
Taking $\Cer>0$ sufficiently large, we obtain that
    \eqan{
    \prob\big({\cal G}^{\sss(n)}_1\cap {\cal G}^{\sss(n)}_2\cap ({\cal G}^{\sss(n)}_3)^c\big)
    &\leq \sum_{i=1}^n \sum_j\prob\big(\Der_i\geq \Cer j^{\tau-1}|D_i=j, {\cal G}^{\sss(n)}_1\cap {\cal G}^{\sss(n)}_2)\prob(D_i=j, {\cal G}^{\sss(n)}_1\cap {\cal G}^{\sss(n)}_2)\nn\\
    &\leq \sum_{i=1}^n {\mathrm e}^{-RCn/i}\leq \vep/3,
    }
where we use the fact that $D_i\geq C^{-1}(n/i)^{1/(\tau-1)}$ since the event
${\cal G}^{\sss(n)}_2$ occurs, and where, in the last step, we use the fact that,
for each $a,C>0$, we can make $R(a,C,\Cer)$ large by taking $\Cer$ sufficiently large.
\qed
\medskip

\noindent Now we are ready to prove the tightness of the FPP
problem. In the statement below, we let ${\cal A}_{\sss
m,K}^{\sss(n)}(i,j)$ be the event that there exists a path of length
$2m$ connecting $i$ and $j$ of weight at most $W/\sqrt{n}$ that
leaves $[K]$, and we write for $k$ fixed,
    \eqn{
    {\cal A}_{\sss m,K}^{\sss(n)}=\bigcup_{i,j\in [k]}{\cal A}_{\sss m,K}^{\sss(n)}(i,j),
    \qquad {\cal A}_{\sss K}^{\sss(n)}=\bigcup_{m=1}^{\infty}{\cal A}_{\sss m,K}^{\sss(n)}.
    }

\begin{Lemma}[Tightness of even shortest-weight paths on the super vertices]
\label{lem-two-paths-tightness}
Fix $k, K\in \Zbold^+$ and $i,j\in [k]$. Then, there exists a $C>0$ such that
    \eqn{
    \prob({\cal A}_{\sss m,K}^{\sss(n)}\cap {\cal G}^{\sss(n)})\leq
    CW K^{-1} {\mathrm e}^{CW\sqrt{\log{K}}}.
    }
\end{Lemma}

\proof
We follow the same line of argument as in the proof of Lemma \ref{lem-two-paths-odd},
but we need to be more careful in estimating the expected number of paths between the super vertices
$i$ and $j$. For $m\geq 2$ and $\veciota=(i_1, \ldots, i_{m-1})$, let $q_{m,n}(\veciota)$ be the expected
number paths with $2m$ edges ($2m$ step paths) such that the position of the path at time $2k$ is equal to $i_k$,
where, by convention, $i_0=i$ and $i_{m}=j$.
Then, similarly as in \eqref{eqn:bmn} but note that now $q_{m,n}(\veciota)$
is an expectation and not a probability, we have that
    \eqn{
    \label{eqn:Amn}
    \prob({\cal A}_{\sss m,K}^{\sss(n)}\cap {\cal G}^{\sss(n)})\leq \frac{CW^{2m}}{(2m)! n^{m}}
    \sum_{\veciota} q_{m,n}(\veciota).
    }
Observe that
    \eqn{
    \lbeq{qmn-bound}
    q_{m,n}(\veciota)\leq \expec[\prod_{s=1}^m \Ner_{i_{s-1}i_s}\indic{{\cal G}^{\sss(n)}}].
    }
%
We further note that, by Lemma \ref{lem-degrees-in-erased}
and on ${\cal G}^{\sss(n)}$,
    \eqn{
    N_{ij}\leq \Der_i\wedge \Der_j\leq \Der_{i\vee j}\leq \Cer n/(i\vee j),
    }
where we abbreviate, for $x,y\in {\mathbb R}$, $x\wedge y=\min\{x,y\}$ and $x\vee y=\max\{x,y\}$.
Thus, by \refeq{qmn-bound}, we arrive at
    \eqn{
    q_{m,n}(\veciota)\leq \prod_{s=1}^m \Cer n/(i_s\vee i_{s-1}),
    }
and hence, after possibly enlarging $\Cer$,
    \eqn{
    \lbeq{AmKn-bd}
    \prob({\cal A}_{\sss m,K}^{\sss(n)}\cap {\cal G}^{\sss(n)})\leq \frac{CW^{2m}}{(2m)! n^{m}}
    \sum_{\veciota} \prod_{s=1}^m \frac{\Cer n}{i_s\vee i_{s-1}}
    =\frac{(\Cer)^m W^{2m}}{(2m)!}
    \sum_{\veciota} \prod_{s=1}^m \frac{1}{i_s\vee i_{s-1}},
    }
where the sum over $\veciota$ is such that there exists at least one $s$ such that
$i_s> K$, because the path is assumed to leave $[K]$. We now bound
\refeq{AmKn-bd}. Let $1\le t\le m$ be such that $i_t=\max_{s=1}^m i_s$, so that $i_t> K$.
Then, using that both $i_s$ and $i_{s-1}$ are smaller than $i_s\vee i_{s-1}$, we can bound
    \eqn{
    \prod_{s=1}^m \frac{1}{i_s\vee i_{s-1}}
    =\Big(\prod_{s=1}^{t-1} \frac{1}{i_s\vee i_{s-1}}\Big)\frac{1}{i_{t-1}\vee i_{t}}
    \frac{1}{i_t\vee i_{t+1}}\Big(\prod_{s=t+2}^{m} \frac{1}{i_s\vee i_{s-1}}\Big)
    \leq \frac{1}{i_t^2} \prod_{s=1}^{t-1} \frac{1}{i_s}
    \prod_{s=t+2}^m \frac{1}{i_{s-1}}.
    }
Thus,
    \eqn{
    \sum_{\veciota} \prod_{s=1}^m \frac{1}{i_s\vee i_{s-1}}
    =\sum_{t=1}^m \sum_{i_t> K} \frac{1}{i_t^2} \sum_{i_1, \ldots, i_{t-1}\leq i_t}\prod_{s=1}^{t-1} \frac{1}{i_s}
    \sum_{i_{t+1}, \ldots, i_{m-1}\leq i_t}\prod_{s=t+2}^m \frac{1}{i_{s-1}}
    =m\sum_{u> K} \frac{1}{u^2} h_u^{m-2},
    }
where
    \eqn{
    h_u=\sum_{v=1}^u \frac{1}{v}.
    }
We arrive at
    \eqn{
    \lbeq{AmnK-bd}
    \prob({\cal A}_{\sss m,K}^{\sss(n)}\cap {\cal G}^{\sss(n)})\leq \frac{(\Cer)^m W^{2m}}{(2m)!}
    m\sum_{u> K} \frac{1}{u^2} h_u^{m-1}.
    }
By Boole's inequality and \refeq{AmnK-bd}, we obtain, after replacing $\Cer$ by $C$, that
    \eqan{
    \prob({\cal A}_{\sss K}^{\sss(n)}\cap {\cal G}^{\sss(n)})
    &\leq \sum_{m=2}^{\infty} \prob({\cal A}_{\sss m,K}^{\sss(n)})
    \leq\sum_{m=2}^{\infty}\frac{C^m W^{2m}}{(2m)!}
    m\sum_{u> K} \frac{1}{u^2} h_u^{m-1}\nn\\
    &\leq W\sum_{u> K} \frac{1}{u^2} \sum_{m=2}^{\infty} h_u^{(2m-1)/2}
    \frac{C^{2m-1} W^{2m-1}}{(2m-1)!}\nn\\
    &\leq CW\sum_{u> K} \frac{1}{u^2} {\mathrm e}^{C W\sqrt{h_u}}
    \leq CW K^{-1} {\mathrm e}^{C W\sqrt{\log{K}}},
    }
where we used that
    \eqan{
    \sum_{u> K} \frac{1}{u^2} {\mathrm e}^{C W\sqrt{h_u}}
    &\le
    \int_{\log K}^\infty {\mathrm e}^{-y+C W\sqrt{y+c}}\,dy\nn\\
    &\leq
    \int_{\log K}^\infty \exp\big\{-y\big(1-\frac{ 2CW}{\sqrt{\log{K}}}\big)\big\}\,dy
    \leq
    K^{-1} {\mathrm e}^{C' W\sqrt{\log{K}}},
    }
for some $c, C'>0.$ This completes the proof of Lemma \ref{lem-two-paths-tightness}.
\qed
\medskip

\noindent
Now we are ready to complete the proof of Proposition \ref{prop-hopcount-conv}:\\
{\it Proof of Proposition \ref{prop-hopcount-conv}.} We write $H_n(i,j)$
and $W_n(i,j)$ for the number of edges and weight of the shortest-weight
path between the super vertices $i,j\in[k]$.

(a) The fact that $\prob(H_n(i,j)\not\in 2\Zbold^+)=o(1)$ for any super vertices $i,j$,
follows immediately from Lemma \ref{lem-two-paths-odd}, which implies that
the even length path between $i$ and $j$ is a two-edge path {\bf whp}.
The tightness of $H_n(i,j)$ follows from part (b),
which we prove next.

(b) By Proposition \ref{prop-weak-con-FPP-erased}, the rescaled weight
$\sqrt{n} W_n(i,j)\leq l_{ij}^{\sss(n)}$ is a tight sequence of
random variables, so that,for $W$ large, it is at most $W$ with probability converging to 1 when
$W\rightarrow \infty$. Fix $\vep>0$ arbitrary.
Then, fix $W>0$ sufficiently large such that
the probability that $\sqrt{n} W_n(i,j)>W$ is at most $\vep/3$, $K>0$
such that $CW K^{-1} {\mathrm e}^{CW\sqrt{\log{K}}}<\vep/3$,
and, use Lemma \ref{lem-degrees-in-erased} to see that we can choose
$a,C, \Cer$ such that $\prob({\cal G}^{\sss(n)})\geq 1- \vep/3$. Then,
by Lemma \ref{lem-two-paths-tightness}, the probability that this
two-edge path leaves $[K]$ is at most $\vep/3+\vep/3+\vep/3=\vep$.
This completes the proof of (b).

(c) The proof that $\prob(H_n\not\in 2\Zbold^+)=o(1)$ follows from (a) since,
for $k$ large, {\bf whp}, $A_1$ and $A_2$ are exclusively attached to
super vertices in $[k]$. The tightness of $H_n$ also follows from
this argument and (a).
\qed

\subsection{The FPP on $\KKer$ is well defined: Proof of Proposition \ref{prop-infi-FPP-Ker-well-defined}}
\label{sec-pf-prop6.4}
In this section, we prove Proposition \ref{prop-infi-FPP-Ker-well-defined}.
For this, we start by investigating $f(P_i,P_j)$ for large $i,j$.
The main result is contained in the following lemma:

\begin{Lemma}[Asymptotics for $f(P_i,P_j)$ for large $i,j$]
\label{lem-asy-fpipj}
Let $\eta$ be a stable random variable with parameter $\tau-1\in (0,1)$. Then, there exists a
constant $c>0$ such that, as $i\wedge j\rightarrow\infty$,
    \eqn{
    \lbeq{fpipj-asy}
    f(P_i,P_j)\leq \frac{c\eta^{1-\tau}}{i\vee j}, \quad a.s.
    }

\end{Lemma}

\proof We note that, by \refeq{EVTexpsrep} and \refeq{PRPDdef} and the strong law of large numbers that,
as $i\rightarrow \infty$,
    \eqn{
    \eta P_i i^{1/(\tau-1)}=(i/\Gamma_i)^{1/(\tau-1)}\convas 1.
    }
Further, by \refeq{fpipj-def},
    \eqn{
    \lbeq{fst-asy}
    f(s,t)=1-\expec[(1-s)^D]-\expec[(1-t)^D]+\expec[(1-s-t)^D]
    \leq 1-\expec[(1-s)^D] \leq cs^{\tau-1},
    }
since, for $\alpha=\tau-1,$ and $D$ in the domain of attraction of an $\alpha$-stable random variable, we have that,
as $u\downarrow 0$,
    \eqn{
    \expec[{\mathrm e}^{-u D}]={\mathrm e}^{-c u^\alpha(1+o(1))}=1-cu^{\alpha}(1+o(1)).
    }
Combing these asymptotics proves \refeq{fpipj-asy}.
\qed
\medskip

\noindent
{\it Proof of Proposition \ref{prop-infi-FPP-Ker-well-defined}.}
Let ${\cal A}_{\sss m,K}$ be the event that there exists a path
of length $m$ and weight at most $W$ connecting $i$ and $j$ and which contains a
vertex in $\Zbold^+\setminus [K]$.
Then, by Boole's inequality and the conditional independence
of the weights $\{l_{ij}\}_{1\leq i\leq j<\infty}$, we obtain that
    \eqn{
    \prob({\cal A}_{\sss m,K})
    \leq \sum_{\veciota} \prob(\sum_{s=1}^m l_{i_{s-1} i_{s}} \leq W),
    }
where, as in the proof of Lemma \ref{lem-two-paths-tightness}, the sum over
$\veciota$ is over $\veciota=(i_1, \ldots, i_{m-1})$, where, by convention, $i_0=i$ and $i_{m}=j$,
and $\max_{s=1}^m i_s\geq K$. Now, by the conditional independence of $\{l_{ij}\}_{1\leq i<j<\infty}$,
    \eqan{
    \prob(\sum_{s=1}^m l_{i_{s-1} i_{s}} \leq W|\{P_i\}_{i\geq 1})
    &=\int_{x_1+\cdots+x_m\leq W} \prod_{s=1}^m f(P_{i_{s-1}},P_{i_s}) x_s {\mathrm e}^{-f(P_{i_{s-1}},P_{i_s}) x_s^2/2} d x_1\cdots dx_m\nn\\
    &\leq   \prod_{s=1}^m f(P_{i_{s-1}},P_{i_s}) \int_{x_1+\cdots+x_m\leq W} x_1\cdots x_m d x_1\cdots dx_m\nn\\
    &=\frac{W^{2m}}{(2m)!} \prod_{s=1}^m f(P_{i_{s-1}},P_{i_s}),
    }
by \cite[4.634]{GraRyz65}. We have that $f(P_i,P_j)\leq c\eta^{1-\tau}(i\vee j)^{-1}, \, a.s.,$
by Lemma \ref{lem-asy-fpipj}. The random variable $\eta$ has a stable
distribution, and is therefore {\bf whp} bounded above by $C$ for some $C>0$
sufficiently large. The arising bound is identical to the bound \refeq{AmKn-bd} derived in
the proof of Lemma \ref{lem-two-paths-tightness}, and we can follow
the proof to obtain \refeq{infi-FPP-bound}.
\qed

\section{Robustness and fragility: Proof of Theorem \ref{thm-rob-frag}}
\label{sec-rob-frag-pf}

We start by proving Theorem \ref{thm-rob-frag}(a), for which we note that whatever
the value of $p\in (0,1)$, {\bf whp}, not all super vertices will be deleted.
The number of undeleted vertices that are connected to a kept super vertex will be
$\Theta_{\sss \prob}(n)$, which proves the claim. In fact, we now argue that a stronger
result holds. We note that the size of the giant component is the same wether we
consider $\GGer_n$ or $\GGor_n$. It is easy to prove the following result:

\begin{Theorem}[Giant component after random attack]
Consider either $\GGor_n$ or $\GGer_n$ and leave each vertex with probability $p$ or delete it with probability $(1-p)$. The resulting graph (of vertices which are left) has a unique giant component $\CC_n(p)$. Further,
with $|\CC_n(p)|$ denoting the number of vertices in $\CC_n(p)$,
    \eqn{
    \lbeq{Var-GC-CM}
    \frac{\expec[|\CC_n(p)|]}{n} \longrightarrow p \expec[1 - (1-p)^{\Der}] = \lambda(p),
    \qquad
    \mbox{Var}\left(\frac{|\CC_n(p)|}{n}\right) \to \beta(p)> 0.
    }
\end{Theorem}

Unlike for other random graph models, \refeq{Var-GC-CM} suggests that $|\CC_n(p)|/n \convd Z_p$, where $Z_p$
is a \emph{non-degenerate} random variable. We shall however not attempt to prove the latter statement here.

\noindent
{\bf Sketch of proof:} Note that we have the identity
    \[
    \frac{\expec[|\CC_n(p)|]}{n} = \prob(1 \in \CC_n(p)),
    \]
where $1$ is a uniformly chosen vertex in $\GGer_n$. For large $n$, the vertex $1$ being in the giant component is {\it essentially} equivalent to the following two conditions:
\\(i) Vertex $1$ is not deleted;  this happens with probability $p$.
\\(ii) Vertex $1$ is attached to $\Der$ super vertices. If one of those super vertices is not deleted,
then the component of this super vertex is of order $n$ and thus has to be the giant component.
Thus at least one of the super vertices to which $1$ is attached should remain undeleted;
conditionally on $\Der$, this happens with probability $1 - (1-p)^{\Der}$.
Combining (i) and (ii) gives the result. A calculation, using similar ideas as in the proof of
Lemma \ref{lem-degrees-in-erased}, suggests that $\lambda(p)=\Theta(p^2)$ when $p\downarrow 0.$
Further, the giant component is unique, since any pair of super vertices which are kept
are connected top each other, and are each connected to $\Theta_{\sss \prob}(n)$ other vertices.

To prove the convergence of the variance, we note that
    \eqn{
    \mbox{Var}(|\CC_n(p)|)=\sum_{i,j}
    \Big[\prob(i,j\in \CC_n(p))-\prob(i\in \CC_n(p))\prob(j\in \CC_n(p))\Big].
    }
Thus,
    \eqn{
    \mbox{Var}(|\CC_n(p)|/n)=\prob(1,2\in \CC_n(p))-\prob(1\in \CC_n(p))\prob(2\in \CC_n(p)),
    }
where $1,2$ are two independent uniform vertices in $[n]$. Now,
    \eqn{
    \prob(1,2\in \CC_n(p))
    =p^2 \prob(1,2\in \CC_n(p)|1,2 \text{ kept})+o(1),
    }
and
    \eqan{
    &\prob(1,2\in \CC_n(p)|1,2 \text{ kept})
    =1-\prob(\{1\not \in \CC_n(p)\}\cup \{2\not \in \CC_n(p)\}|1,2 \text{ kept})\\
    &\quad=1-\prob(1\not \in \CC_n(p)|1,2 \text{ kept})-\prob(2\not \in \CC_n(p)|1,2 \text{ kept})+\prob(1,2 \not \in \CC_n(p)|1,2 \text{ kept})\nn\\
    &\quad=1-\prob(1\not \in \CC_n(p)|1 \text{ kept})-\prob(2\not \in \CC_n(p)|2 \text{ kept})+\prob(1,2 \not \in \CC_n(p)|1,2 \text{ kept})+o(1),\nn
    }
so that
    \eqn{
    \mbox{Var}(|\CC_n(p)|/n)=
    p^2\prob(1,2 \not \in \CC_n(p)|1,2 \text{ kept})-p^2\prob(1\not \in \CC_n(p)|1
    \text{ kept})\prob(2\not \in \CC_n(p)|2 \text{ kept})+o(1).
    }
Then, we compute that
    \eqn{
    \prob(1\not \in \CC_n(p)|1 \text{ kept})=\expec[(1-p)^{\Der_1}],
    }
while
    \eqn{
    \prob(1,2 \not \in \CC_n(p)|1,2 \text{ kept})
    =\expec[(1-p)^{\Der_1+\Der_2-\Ner_{12}}],
    }
where $\Der_1, \Der_2$ are conditionally independent given $\{P_i\}_{i\geq 1}$, and
$\Ner_{12}$ denotes the number of joint neighbors of $1$ and $2$, and we
use that the total number of super vertices to which 1 and 2 are connected
is equal to $\Der_1+\Der_2-\Ner_{12}$. As a result,
    \eqn{
    \mbox{Var}(|\CC_n(p)|/n)=p^2\Big(\expec[(1-p)^{\Der_1+\Der_2-\Ner_{12}}]
    -\expec[(1-p)^{\Der_1}]\expec[(1-p)^{\Der_2}]\Big)+o(1),
    }
which identifies
    \eqn{
    \beta(p)=p^2\Big(\expec[(1-p)^{\Der_1+\Der_2-\Ner_{12}}]
    -\expec[(1-p)^{\Der_1}]\expec[(1-p)^{\Der_2}]\Big).
    }
To see that $\beta(p)>0$, we note that $\Ner_{12}>0$ with positive probability, so that
    \eqn{
    \expec[(1-p)^{\Der_1+\Der_2-\Ner_{12}}]
    >\expec[(1-p)^{\Der_1+\Der_2}]=\expec\Big[\expec[(1-p)^{\Der_1+\Der_2}\mid \{P_i\}_{i\geq 1}]
    \Big]=\expec\Big[\expec\big[(1-p)^{\Der_1}\mid \{P_i\}_{i\geq 1}\big]^2\Big],
    \nn
    }
by the conditional independence of $\Der_1$ and $\Der_2$. Thus, $\beta(p)>0$
by the Cauchy-Schwarz inequality, as claimed.
\qed

To prove Theorem \ref{thm-rob-frag}(b), we again use that a uniform vertex is, {\bf whp},
only connected to a super vertex. Thus, there exists $K_\eps$ such that
by deleting the $K_\eps$ vertices with largest degree,
we shall isolate $A_1$ with probability at least $\vep$. This proves
\refeq{disconnection-prob}.
\qed

\section{Conclusion}
\label{sec-conc}
We conclude with a 
discussion about various extensions of
the above results together with some further results without proof. Throughout the
discussion we shall use $\GG_n$ to denote either of $\GGor_n$ and $\GGer_n$,
where the choice depends on the context under consideration.
\\(a) {\bf Load distribution:} Understanding how random disorder changes
the geometry of the network is crucial for understanding asymptotics of
more complicated constructs such as the load distribution. More precisely,
for any pair of vertices $i,j\in \GG_n$,  let $\bpi(i,j)$ denote the
minimal weight path between the two vertices. For any vertex $v\in \GG_n$,
the {\it load} on the vertex is defined as
    \[
    L_n(v) = \sum_{i\neq j} \ind_{\{v \in \bpi(i,j)\}}.
    \]
For any fixed $x$, define the function $G_n(x)$ as
    \[
    G_n(x)  = \#\{v: L_n(v) > x \}.
    \]
Understanding such functions is paramount to understanding the flow carrying
properties of the network and are essential for the study of {\it betweenness centrality}
of vertices in a network. For example in social networks, such measures are used to rate the
relative importance of various individuals in the network, while in data networks such as
the World-Wide Web, such measures are used to rank the relative importance of web pages.
An actual theoretical analysis of such questions is important but seems difficult in many
relevant situations. It would be of interest to find asymptotics of such functions
in terms of the infinite objects $\KKor$ and $\KKer$ constructed in this paper.
See also \cite{SD07} for an analysis of such questions in the mean-field setting.
\\(b) {\bf Universality for edge weights:} In this study, to avoid technical
complications we assumed that each edge weight in $\GG_n$ has an exponential
distribution.   One natural question is how far do these results depend on
this assumption. It is well known in probabilistic combinatorial optimization
that in a  wide variety of contexts, when considering problems such as those
in this paper, the actual distribution of the edge weights is
not that important, what is important is the value of the density at $0$.
More precisely, consider $\GGer_n$ (i.e., the erased CM)
where each edge is given an i.i.d.\ edge weight having a continuous
distribution  with density $g$ and let $g(0) = \zeta\ch{\in (0,\infty)}$. Similar to
$\KKer$ defined in Section \ref{sec-erased}, define $\KKer(\zeta)$ to be the
infinite graph on the vertex set $\Zbold^+$ where each edge $l_{ij}$ has the distribution
    \eqn{
    \lbeq{lij-gen}
    \prob\left(l_{ij}> x\right) = \exp\left(-f(P_i, P_j)\zeta^2 x^2/2\right).
    }
Equation \refeq{lij-gen} can be proved along similar lines as
in the proof of Lemma \ref{lem-min-Gamma}, and we leave this to the reader.

Let $\Ier$ be as defined in Section \ref{sec-erased}.  Then we have the following
modification of Theorem \ref{thm-erased} which can be proved along the same lines:

\begin{Theorem}[Extension to other densities]
\label{thm-erased-ext}
Theorem \ref{thm-erased} continues to hold with the modification that the quantities $\Wer_{ij}, \Her_{ij}$
arising in the limits are replaced by the corresponding quantities in $\KKer(\zeta)$ instead of $\KKer$, $\Ver_i$ is
distributed as the minimum of $\Der$ random variables having density $g$,  while the distributions of  $\Ier$
and $\Jer$ remain unchanged.
\end{Theorem}

\noindent
\ch{A more challenging extension would be to densities for which either $g(0)=0$, or for which
$\lim_{x\downarrow 0} g(x)=\infty$. In this case, we believe the behavior to be entirely
different from the one in Theorems \ref{thm-erased} and \ref{thm-erased-ext}, and
it would be of interest to investigate whether a similar limiting FPP process arises.}

\paragraph{Acknowledgments.}
The research of SB is supported by N.S.F. Grant DMS 0704159, NSERC and PIMS Canada.
SB would like to thank the hospitality of Eurandom where much of this work was done.
The work of RvdH was supported in part by Netherlands Organisation for Scientific
Research (NWO).

\bibliographystyle{plain}
\bibliography{bib-sparse}

\def\cprime{$'$}
\begin{thebibliography}{10}

\bibitem{ba-rob}
R.~Albert, H.~Jeong, and A.L. Barab{\'a}si.
\newblock {Error and attack tolerance of complex networks}.
\newblock {\em Nature(London)}, {\bf 406}(6794):378--382, (2000).

\bibitem{SD07}
D.J. Aldous and S.~Bhamidi.
\newblock {Flows through random networks}.
\newblock {\em To appear in Random Structures and Algorithms}, (2009).

\bibitem{BenCan78}
E.A. Bender and E.R. Canfield.
\newblock The asymptotic number of labelled graphs with a given degree
  sequences.
\newblock {\em Journal of Combinatorial Theory (A)}, {\bf 24}:296--307, (1978).

\bibitem{vcg-random-shanky}
S.~Bhamidi.
\newblock First passage percolation on locally tree like networks {I}: Dense
  random graphs.
\newblock {\em Journal of Mathematical Phyiscs}, {\bf 49}:125218, (2008).

\bibitem{BhaHofHoo09b}
S.~Bhamidi, R.~van~der Hofstad, and G.~Hooghiemstra.
\newblock First passage percolation on sparse random graphs with finite mean
  degrees.
\newblock Available from {\tt http://arxiv.org/abs/0903.5136}, Preprint (2009).

\bibitem{Boll01}
B.~Bollob{\'a}s.
\newblock {\em Random graphs}, volume~{\bf 73} of {\em Cambridge Studies in
  Advanced Mathematics}.
\newblock Cambridge University Press, Cambridge, second edition, (2001).

\bibitem{boll-rior-perc}
B.~Bollob{\'a}s and O.~Riordan.
\newblock {Robustness and Vulnerability of Scale-Free Random Graphs}.
\newblock {\em Internet Mathematics}, {\bf 1}(1):1--35, (2004).

\bibitem{DeiEskHofHoo09}
M.~Deijfen, H.~van~den Esker, R.~van~der Hofstad, and G.~Hooghiemstra.
\newblock A preferential attachment model with random initial degrees.
\newblock {\em Arkiv f{\"o}r Matematik}, {\bf 47}:41--72, (2009).

\bibitem{durrett-lnp}
R.~Durrett.
\newblock {Lecture notes on particle systems and percolation}.
\newblock {\em Pacific Grove, CA}, 1988.

\bibitem{hofs2}
H.~van~den Esker, R.~van~der Hofstad, G.~Hooghiemstra, and D.~Znamenski.
\newblock Distances in random graphs with infinite mean degrees.
\newblock {\em Extremes}, {\bf 8}(3):111--141, (2005).

\bibitem{GraRyz65}
I.~S. Gradshteyn and I.~M. Ryzhik.
\newblock {\em Table of integrals, series, and products}.
\newblock Fourth edition prepared by Ju. V. Geronimus and M. Ju. Ce\u\i tlin.
  Translated from the Russian by Scripta Technica, Inc. Translation edited by
  Alan Jeffrey. Academic Press, New York, (1965).

\bibitem{hamm-welsh}
J.M. Hammersley and D.J.A. Welsh.
\newblock {First-passage percolation, sub-additive process, stochastic network
  and generalized renewal theory}.
\newblock {\em Bernoulli, 1713: Bayes, 1763; Laplace, 1813. Anniversary
  Volume}, (1965).

\bibitem{hofs-erdos-fpp}
R.~van~der Hofstad, G.~Hooghiemstra, and P.~Van~Mieghem.
\newblock First-passage percolation on the random graph.
\newblock {\em Probab. Engrg. Inform. Sci.}, {\bf 15}(2):225--237, (2001).

\bibitem{hofs-flood}
R.~van~der Hofstad, G.~Hooghiemstra, and P.~Van~Mieghem.
\newblock The flooding time in random graphs.
\newblock {\em Extremes}, {\bf 5}(2):111--129 (2003), (2002).

\bibitem{hofs3}
R.~van~der Hofstad, G.~Hooghiemstra, and P.~Van~Mieghem.
\newblock Distances in random graphs with finite variance degrees.
\newblock {\em Random Structures Algorithms}, {\bf 27}(1):76--123, (2005).

\bibitem{hofs1}
R.~van~der Hofstad, G.~Hooghiemstra, and D.~Znamenski.
\newblock Distances in random graphs with finite mean and infinite variance
  degrees.
\newblock {\em Electron. J. Probab.}, {\bf 12}:no. 25, 703--766 (electronic),
  (2007).

\bibitem{howard}
C.D. Howard.
\newblock {Models of first-passage percolation}.
\newblock {\em Probability on Discrete Structures}, pages 125--173, (2004).

\bibitem{janson123}
S.~Janson.
\newblock One, two and three times {$\log n/n$} for paths in a complete graph
  with random weights.
\newblock {\em Combin. Probab. Comput.}, {\bf 8}(4):347--361, (1999).
\newblock Random graphs and combinatorial structures (Oberwolfach, 1997).

\bibitem{Jans02}
S.~Janson.
\newblock On concentration of probability.
\newblock {\em Contemporary Combinatorics, ed. B. Bollob{\'a}s, Bolyai Soc.
  Math. Stud.}, {\bf 10}:289--301, (2002).
\newblock J{\'a}nos Bolyai Mathematical Society, Budapest.

\bibitem{Kall02}
O.~Kallenberg.
\newblock {\em Foundations of modern probability}.
\newblock Probability and its Applications (New York). Springer-Verlag, New
  York, second edition, (2002).

\bibitem{babak-newman2}
L.A. Meyers, M.E.J. Newman, and B.~Pourbohloul.
\newblock {Predicting epidemics on directed contact networks}.
\newblock {\em Journal of Theoretical Biology}, {\bf 240}(3):400--418, (2006).

\bibitem{babak-newman1}
L.A. Meyers, B.~Pourbohloul, M.E.J. Newman, D.M. Skowronski, and R.C. Brunham.
\newblock {Network theory and SARS: predicting outbreak diversity}.
\newblock {\em Journal of Theoretical Biology}, {\bf 232}(1):71--81, (2005).

\bibitem{MolRee95}
M.~Molloy and B.~Reed.
\newblock A critical point for random graphs with a given degree sequence.
\newblock In {\em Proceedings of the Sixth International Seminar on Random
  Graphs and Probabilistic Methods in Combinatorics and Computer Science,
  ``Random Graphs '93'' (Pozna\'n, 1993)}, volume~{\bf 6}, pages 161--179,
  (1995).

\bibitem{newman}
M.E.J. Newman.
\newblock {The structure and function of complex networks}.
\newblock {\em Arxiv preprint cond-mat/0303516}, (2003).

\bibitem{NorRei04}
I.~Norros and H.~Reittu.
\newblock On the power-law random graph model of massive data networks.
\newblock {\em Performance Evaluation}, {\bf 55}:3--23, (2004).

\bibitem{PitYor97}
J.~Pitman and M.~Yor.
\newblock The two-parameter {P}oisson-{D}irichlet distribution derived from a
  stable subordinator.
\newblock {\em Ann. Probab.}, {\bf 25}(2):855--900, (1997).

\bibitem{wastlund}
J.~W{\"a}stlund.
\newblock {Random assignment and shortest path problems}.
\newblock In {\em Fourth Colloquium on Mathematics and Computer Science
  Algorithms, Trees, Combinatorics and Probabilities", DMTCS Proceedings},
  pages 31--38. (2006).
\newblock Available from {\tt
  http://www.dmtcs.org/dmtcs-ojs/index.php/proceedings/issue/view/84/showToc}.

\end{thebibliography}

\end{document}